\documentclass{amsart}

\usepackage{mathrsfs}
\usepackage{amsmath}
\usepackage{amssymb}
\usepackage{graphicx}
\usepackage{color}
\usepackage{caption,subcaption}
\usepackage{babel}
\usepackage{comment}
\usepackage{algorithm}
\usepackage{algpseudocode}
\usepackage{setspace}

\usepackage{grffile}
\usepackage{enumerate}

\usepackage{hyperref}
\hypersetup{colorlinks=true,
            linkcolor=black,
            anchorcolor=black,
            citecolor=black}
            
\newtheorem{theorem}{Theorem}[section]

\theoremstyle{definition}

\theoremstyle{remark}
\newtheorem{remark}[theorem]{Remark}

\numberwithin{equation}{section}

\renewcommand{\a}{\alpha}

\newcommand{\R}{\mathbb{R}}

\newcommand{\lam}{\lambda}

\newcommand{\bx}{\mathbf{x}}

\newcommand{\bu}{\mathbf{u}}

\newcommand{\be}{\mathbf{e}}
\newcommand{\sig}{\sigma}
\newcommand{\gam}{\gamma}
\newcommand{\bv}{\mathbf{v}}
\renewcommand{\th}{\theta}
\newcommand{\tg}{\widetilde{g}}
\newcommand{\tG}{\widetilde{G}}
\newcommand{\tu}{\widetilde{u}}
\newcommand{\tU}{\widetilde{U}}

\newcommand{\lan}{\langle}
\newcommand{\ran}{\rangle}
\newcommand{\eps}{\varepsilon}
\newcommand{\pressure}{\mathbb{P}}

\newcommand{\network}{\mathcal{N}}
\newcommand{\bm}{\mathbf{\mathfrak{m}}}

\newcommand{\hm}{\widehat{m}}
\newcommand{\bg}{\mathbf{g}}
\newcommand{\bE}{\mathbb{E}}
\newcommand{\bB}{\mathbb{B}}

\newcommand{\bR}{\mathbf{R}}
\newcommand{\bn}{\mathbf{n}}
\newcommand{\bq}{\mathbf{q}}
\newcommand{\Dx}{\Delta x}
\newcommand{\Dt}{\Delta t}

\newcommand{\rd}{\mathrm{d}}

\newcommand{\magenta}{\color{magenta}}
\newcommand{\nc}{\normalcolor}
\renewcommand{\tilde}[1]{\widetilde{#1}}

\DeclareMathOperator*{\trace}{trace}
\DeclareMathOperator*{\image}{Im}
\DeclareMathOperator*{\kernel}{ker}

\numberwithin{equation}{section}
\numberwithin{figure}{section}
\numberwithin{table}{section}

\newtheorem{THEOREM}{Theorem}[section]

\newtheorem{proposition}[THEOREM]{Proposition}

\usepackage{xparse}
\makeatletter
\RenewDocumentCommand{\title}{om}{%
   \IfNoValueTF{#1}
     {\gdef\shorttitle{Continuous data assimilation for hydrodynamics}}
     {\gdef\shorttitle{#1}}%
   \gdef\@title{#2}%
}
\makeatother

\begin{document}

\title{Continuous data assimilation for hydrodynamics: consistent discretization and application to moment recovery}

\author{Jingcheng Lu}
\address{School of Mathematics, Univerisity of Minnesota--Twin Cities, Minneapolis, Minnesota 55455}
\email{lu000688@umn.edu}

\author{Kunlun Qi}
\address{School of Mathematics, Univerisity of Minnesota--Twin Cities, Minneapolis, Minnesota 55455}
\email{kqi@umn.edu, kunlunqi.math@gmail.com}

\author{Li Wang}
\address{School of Mathematics, Univerisity of Minnesota--Twin Cities, Minneapolis, Minnesota 55455}
\email{liwang@umn.edu}

\author{Jeff Calder}
\address{School of Mathematics, Univerisity of Minnesota--Twin Cities, Minneapolis, Minnesota 55455}
\email{jwcalder@umn.edu}
\thanks{JL and JC were supported by NSF-CCF:2212318, and JC was additionally supported by an Albert and Dorothy Marden Professorship and a McKnight Presidential Fellowship. LW is partially supported by NSF grant DMS-1846854 and UMN DSI-SSG-4886888864. KQ is partially supported by direct funding from School of Mathematics at the University of Minnesota--Twin Cities.}


\date{}

\emergencystretch 3em


\begin{abstract}
Motivated by the challenge of moment recovery in hydrodynamic approximation in kinetic theory, we propose a data-driven approach for the hydrodynamic models. Inspired by continuous data assimilation, our method introduces a relaxation-based nudging system coupled with a novel discretization technique. This approach facilitates the simultaneous recovery of both the force term and a high-resolution solution from sparsely observed data. To address potential numerical artifacts, we use kernel regression to fit the observed data. We also analyze the convergence of the proposed nudging system under both full and partial data scenarios. When applied to moment systems, the source term involves the derivative of higher-order moments, and our approach serves as a crucial step for data preparation in machine-learning based moment closure models. Multiple numerical experiments demonstrate the effectiveness of our algorithm, and we discuss its potential extension to high-dimensional systems.
\end{abstract}

\maketitle

\tableofcontents

\section{Introduction}

When modeling complicated physics processes, an accurate description can be provided at the mesoscopic scale using kinetic-type governing equations:
\begin{equation}\label{eq:vlasov}
    \frac{\partial}{\partial t}f+\bv\cdot\nabla_{\bx} f+F \cdot\nabla_{\bv} f = \mathcal Q(f,f), \quad (\bx, \bv,t) \in\R^3\times\R^3\times\R_+\,.
\end{equation}
Here $f(\bx,\bv,t)$ is the distribution function that describes particles within a physical system, where $\bx$ is the position,  $\bv$ is the velocity and $t$ is the time. $F$ is the external force, which, for example, can take the form $F = q(\bE+\bv\times \bB)$ when modeling charged particles in plasma. Here 
$q$ is the charge per unit mass of the particles, while $\bE$ and $\bB$ denote the strengths of the electric and magnetic fields, respectively. $\mathcal Q(f,f)$  is the collision operator that encapsulates the interactions between particles.

While significant progress has been made in developing computational methods for solving \eqref{eq:vlasov}, such as finite difference or finite element methods, in industrial applications like fusion energy simulations, these grid-based methods are still cursed by dimensionality.
To reduce the computational cost while maintaining important physical features, an attractive approach is to simulate the dynamics at the hydrodynamic scale, which is the primary motivation behind the moment method. To explain, by taking the zeroth and the first $\bv-$moments of the equation \eqref{eq:vlasov}, one obtains the following hydrodynamic system:
\begin{equation}\label{eq:magnetohydro}
\left\{\begin{array}{l}\displaystyle
     \frac{\partial}{\partial t}\rho+\nabla\cdot(\rho\bu) = 0\,,\\
       \\
       \displaystyle
      \frac{\partial}{\partial t}\rho\bu+\nabla\cdot(\rho\bu\otimes\bu) = -\nabla\cdot\pressure+\rho q(\bE+\bu\times\bB)\,,
\end{array}
\right.
\end{equation}
where the density $\rho$, the momentum $\rho\bu$, and the pressure tensor $\pressure$ are given by 
\[
\begin{split}
&\rho(\bx,t) := \int f(\bx,\bv,t) \rd\bv, \\
&\rho\bu(\bx,t) := \int \bv f(\bx,\bv,t) \rd\bv,\\
& \pressure(\bx,t) := \int (\bv-\bu(\bx,t))\otimes  (\bv-\bu(\bx,t))f(\bx,\bv,t)\rd\bv.
\end{split}
\]
Clearly, the hydrodynamic system \eqref{eq:magnetohydro} resides in a 4-dimensional space: 3 for $\bx$ and 1 for $t$, making it more computationally friendly than the original 7-dimensional kinetic equation \eqref{eq:vlasov}. However, this reduction comes at a price--the system is not closed. Specifically, $\mathbb P$ depends on the second moments of $f$, which cannot be written solely as a function of the zero-th and first moment of $f$, i.e., $\rho$ and $\rho \bu$. One way to close this system is to assume $f$ is near the Maxwellian, i.e., $f\propto \rho \exp{-\frac{|\bv-u|^2}{2 T}} $ . However, this approximation holds true only in certain regimes where collisional effects dominate \cite{bardos1993fluid}. In general, one can derive a larger moment system that evolves more moments, aiming to approximate the original kinetic equation more closely, but the closure problem persists \cite{cai2010numerical, cai2012numerical}.


In this paper, we aim to address this problem with the help of observational data. In particular, we reformulate this closure problem into a force-identification problem: 
\begin{equation}\label{eq:hydro}
   \left\{
   \begin{array}{l}\displaystyle
   \frac{\partial}{\partial t}U(\bx,t)+\nabla\cdot F(U(\bx,t)) = G(\bx,t)+S(U,\bx,t)\,,\\
   \\
   \displaystyle
   U(\bx,0) = U_0(\bx)\,,
   \end{array}
   \right. \quad \bx \in\Omega\subseteq \R^d, \ \ t\geq 0.
\end{equation}
Here $U(\bx,t) = (U^0(\bx,t),U^1(\bx,t),\cdots,U^N(\bx,t))^{\top}\in\R^N$ are the solution states. The flux, $F(U)$, and the source term, $S(U,\bx,t)$, are assumed to be known from the physics, whereas the forcing term $G(\bx,t)$ is to be determined. In the context of moment recovery, the force function takes the divergence form $G = -\nabla\cdot\mathbb{M}$, where the tensor $\mathbb M\in \R^{d\times N}$ contains the information of the moment one order higher than $U^N$. Indeed, the system \eqref{eq:magnetohydro} can be written in the form of \eqref{eq:hydro} with
\[
U = \left[\begin{array}{c}
     \rho \\
     \rho\bu
\end{array}\right], \quad
F = \left[\begin{array}{c}
     \rho \bu \\
     \rho\bu\otimes\bu
\end{array}\right], \quad
G = -\nabla\cdot\left[\begin{array}{c}
     \mathbf 0 \\
     \pressure
\end{array}\right],\quad
S = \left[\begin{array}{c}
     \mathbf 0 \\
     \rho q (\bE+\bu\otimes\bB)
\end{array}\right].
\]

If the solution state $U(\bx,t)$ and its gradients are known continuously in space and time, then the forcing term can be computed with $\displaystyle G = \frac{\partial U}{\partial t}+\nabla\cdot F-S$. However, in many real-world applications such as weather prediction, one may only have access to a limited amount of true data at sparse and/or non-uniform observation points, $\bx^{ob}_{j}$, $j=1,2,\cdots,N_{ob}$, which may not be sufficient to recover the solution gradients with the desired accuracy. Additionally, due to the limited observations, the governing equation \eqref{eq:hydro} may not be initialized accurately. This limitation can sometimes lead to a severe accumulation of errors when making future predictions. All of these challenges underscore the pressing task:

\medskip
{\it Given only sparse observations $\{U(\bx^{ob}_j,t)\}$, can we simultaneously recover $G$ and obtain an enhanced resolution of the  solution $U(\bx, t)$ in \eqref{eq:hydro}? }
\medskip


To fulfill this task, we propose the following nudged system:
\begin{subequations}\label{eq:nudge}
\begin{equation}\label{eq:nudge eqn}
 \left\{
   \begin{array}{l}\displaystyle
   \frac{\partial}{\partial t}V(\bx,t)+\nabla\cdot F(V(\bx,t)) = \tG(\bx,t)+S(V,\bx,t)+\mu(I_h(U)-I_h(V)) \,,\\
   \\
   \displaystyle
   V(\bx,0) = V_0(\bx)\,.
   \end{array}
   \right.
\end{equation}
Here $I_h(\cdot)$ is the \emph{interpolation operator}, which is a bounded linear operator that interpolates the observed data into a smooth function, i.e.,  $ I_h: \{(\bx^{ob}_j, U(\bx^{ob}_j))\}  \mapsto I_h(U)(\bx)\in C^{\infty}(\Omega)$. $h$ represents the mesh size of observation grids $\{\bx^{ob}_j\}$. 
The approximate force $\tG$ is constructed with
\begin{equation}\label{eq:nudge force}
\tG = \frac{\partial}{\partial t} I_h(U)+I_h(\nabla\cdot F(\tU)-S(\tU)), \quad \tU = I_h(U)+(I-I_h)(V).
\end{equation}
\end{subequations}
In the particular case where $I_h$ is a projection, we have $I_h^2 = I_h$. Note that $\partial_t I_h = I_h \partial_t$ and $I_h(\tU) = I_h(U)$, \eqref{eq:nudge force} can also be written as 
\begin{equation*}
\tG = I_h(\frac{\partial }{\partial t}\tU+\nabla\cdot F(\tU)-S(\tU)).
\end{equation*}

Compared with the exact dynamics \eqref{eq:hydro}, the nudged equation \eqref{eq:nudge eqn} includes an additional relaxation term, $\mu(I_h(U)-I_h(V))$, which serves as a {\it feedback control} from the observations $I_h(U)$. The user-tuning parameter  $\mu>0$  is known as the \emph{nudging coefficient}. With an appropriate choice of $I_h$, the approximate solutions $V$ and $\tG$ are expected to be time-asymptotically accurate, namely, $\displaystyle V\overset{t\rightarrow \infty}{\longrightarrow} U$ and $\displaystyle \tG\overset{t\rightarrow \infty}{\longrightarrow} G$ in certain norms. 

We note that the concept of introducing feedback control into the system \eqref{eq:nudge eqn} is inspired by the continuous data assimilation methods discussed in \cite{martinez2022reconstruction, farhat2024identifying}, which aim to recover the unknown body force in the incompressible Navier-Stokes equations from low-mode observations. These methods trace back to the earlier work by Azouani, Olson and Titi \cite{azouani2014continuous, gesho2016computational}, which focuses on velocity field recovery for the two-dimensional Navier-Stokes equations. Extension to the continuous data assimilation for dynamics with
a set of unidentified parameters can be found in \cite{clark2018inferring,carlson2020parameter,carlson2021dynamically,martinez2022convergence,pachev2022concurrent}. Investigations have also been carried out for state recovery in various of hydrodynamic settings, including Rayleigh-Benard convection \cite{farhat2015continuous,farhat2016data,altaf2017downscaling,farhat2020data}, geophysical fluids \cite{jolly2017data,albanez2018continuous,pei2018continuous,jolly2019continuous}, and dispersive equations \cite{jolly2015determining,jolly2017determining}.

A major advantage of the Azouani-Olson-Titi nudging framework is that it allows simultaneous reconstruction of both the solution state $U$ and the unknown forcing term $G$. This capability is particularly important in our context, as $G$ encapsulates the information of truncated higher-order moments, which are the critical components for accurate system recovery. In comparison, other conventional data assimilation methods, such as the Kalman filter \cite{kalman1960new,welch1995introduction}, particle filter \cite{kitagawa1996monte,carpenter1999improved}, and variational assimilation \cite{le1986variational,vzupanski1995four,fischer2005overview}, typically assume fully known model dynamics and are thus primarily designed for state estimation. These methods do not inherently support the recovery of unknown or partially known dynamics, and it remains open to explore if they can enable more complicated data recovery tasks with affordable modifications.


To ensure the convergence of \eqref{eq:nudge}, we impose the following conditions on $I_h$:
\begin{subequations}\label{eq:Ih conditions}
\begin{equation}\label{eq:Ih cond 2}
||I_h(\partial^\a\phi)||_2 \leq c_m h^{-m}||\phi||_2, 
\end{equation}
\begin{equation}\label{eq:Ih cond 3}
||\phi-I_h(\phi)||_2 \leq C_m h^m ||\phi||_{H^m}, 
\end{equation}
\end{subequations}
for any multi-index $|\a| = m$ with $m\geq 0$. The constants $c_m, C_m>0$ are independent of $h$. For example, the spectral projection onto Fourier modes with wave numbers $|k|\leq 1/h$ would satisfy the above conditions. In this special case, $I_h$ serves as a low-pass filter that removes the high-frequency fluctuations, and the intermediate solution $\tU$ is obtained with the low-mode observations, $I_h(U)$, corrected by the high-mode discrepancy, $(I-I_h)(V)$. 

\begin{remark}
In the particular case where the function $\phi$ has weak oscillations such that
\[
\sup_{t\geq 0}\frac{||\nabla\phi||_2}{||\phi||_2}\leq C<\infty,
\]
the error bound \eqref{eq:Ih cond 2} can be modified to
\begin{equation}\label{eq: Ih weak oscillation bound}
||I_h(\nabla\phi)||_2\leq c_0||\nabla\phi||_2\leq c_0 C||\phi||_2.
\end{equation}
The $h^{-1}$ dependence in the error bound for $m = 1$ is then removed.

\end{remark}

Upon the recovery of  (the divergence of) the highest moment, an approximate closure relation may be obtained with machine learning techniques. In particular, consider the moment systems of \eqref{eq:vlasov}: 
\begin{equation} \label{eq:rte moment system0}
\begin{split}
    &\frac{\partial \bm_{N}}{\partial t}+A_N\frac{\partial\bm_{N}}{\partial x} = \bg_{N+1}+S_N\bm_{N}, \\
    &\bg_{N+1} = -\frac{N+1}{2N+1}\partial_x m_{N+1} \be_{N+1},
\end{split}
\end{equation}
where $\bm_{N} = (m_0,m_1,\cdots,m_N)^\top$, $\be_{N+1} = (0,\cdots,0,1)^\top$ and 
\[
m_k(x,t) = \frac{1}{2}\int^1_{-1}f(x,v,t)P_k(v) \rd v\,,
\]
with $P_k$ being the $k-$th order polynomial. $A_N$ and $S_N$ are defined in \eqref{eqn:A}
and \eqref{eqn:S}.
Given observations of the first few moments, e.g., $m_1, \cdots, m_l$,  the proposed method allows for the recovery of the higher moments $m_{l+1}, \cdots, m_{N+1}$. By using machine learning techniques, such as neural networks, to approximate the relationship between the gradient of the highest moments and those of the lower-order moments, one can close the system described by \eqref{eq:rte moment system0}. For a simple kinetic equation, such a relation has been learned in a series of works by Huang et. al. \cite{huang2022machineI,huang2023machineII,huang2023machineIII} by running an expensive kinetic simulation to generate data (i.e., $f$), computing its moments, and then training the neural network on these moments. In contrast, our proposed method eliminates the need to compute the full kinetic system, instead providing a cost-efficient mechanism for generating data for the moment closure task.


The rest of the paper is organized as follows. In Section \ref{sec:approach}, we first introduce the main contents of our proposed method, including the convergence analysis of the nudge system, a novel numerical discretization and a kernel regression approach for interpolation. Extensive numerical experiments are presented for the validation of the effectiveness of our algorithm in Section \ref{sec:numerical1}. Furthermore, we apply the proposed approach to the moment system in kinetic theory, aiming to recover higher-order moments  with incomplete observations in Section \ref{sec:moment-recovery}. Finally, in Section \ref{sec:extension}, we discuss the potential strategies and challenges in extending the method to higher-dimensional systems.

\section{Our approach: convergence analysis and numerical discretization}
\label{sec:approach}

This section presents the main approach of this paper, which centers around \eqref{eq:nudge eqn}. In the next subsection, we first analyze the convergence behavior of \eqref{eq:nudge eqn}
under both fully observed data (i.e., when \eqref{eq:observ conditions} is satisfied) and partial data scenarios (see Proposition~\ref{prop:error for sparse observ}). Subsection 2.2 will then propose an efficient discretization of  \eqref{eq:nudge eqn}. The choice of the interpolation operator $I_h$ is then discussed in section 2.3.

\subsection{Convergence analysis}

\begin{proposition}\label{prop:state recovery}\textbf{\emph{(Convergence of the state)}}
Let $(U,G)$ be a strong solution of the equation 
\begin{equation}\label{eq:linear eqn}
\partial_t U +\partial_x U = G+S(U,x,t), 
\end{equation}
associated with Lipschitz continuous source term,
\[
|S(U,x,t)-S(V,x,t)|\leq L_s|U-V|, \quad L_s>0.
\]
Let $(V,\tG)$ be a strong solution of the nudged system 
\begin{equation}\label{eq:linear nudge}
\left\{\begin{array}{l}
\partial_t V+\partial_x V = \tG+S(V)+\mu(I_h(U)-I_h(V))\,,\\
\\
\tG = \partial_t I_h(U)+I_h(\partial_x\tU-S(\tU)), \quad \tU = I_h(U)+(I-I_h)(V)\,,
\end{array}\right.
\end{equation}
where $I_h$ satisfies conditions \eqref{eq:Ih conditions}. Assume that the data are sufficiently observed, meaning that 
\begin{subequations}\label{eq:observ conditions}
\begin{equation}\label{eq:observ cond1}
||W(t)-I_h\big(W(t)\big)||_2 \leq C_h \, e^{-\eps t}, \quad W(t) = U(t)-V(t),
\end{equation}
\begin{equation}\label{eq:observ cond2}
I_h\big(G(t)\big) = G(t), \\
\end{equation}
\end{subequations}
where $C_h$, $\eps>0$ are positive constants. Then the nudged solution $V(x,t)$ recovers the true data $U(x,t)$ at an exponential rate in time,
\begin{equation}\label{eq:state exp decay}
||U(\cdot, t)-V(\cdot, t)||_2 \leq ||U(\cdot, 0)-V(\cdot, 0)||_2 \, e^{-(\mu-L_s)t}+C_{\mu,h} \, e^{-\eps t},
\end{equation}
with $C_{\mu,h} = (c_1h^{-1}+C_0 L_s+\mu)C_h/(\mu-L_s-\eps)$.
\end{proposition}
\begin{remark}\label{rem:mu effect}
From \eqref{eq:state exp decay}, it appears that a larger $\mu$ results in faster convergence. However, in practice, a larger $\mu$ imposes stricter constraints on the time step due to stability considerations.
\end{remark}
\emph{Proof.} Let $W = U-V$, the difference of \eqref{eq:linear eqn} and \eqref{eq:linear nudge} yields
\[
\frac{\partial W}{\partial t}+\frac{\partial W}{\partial x} = G-\tG+S(U)-S(V)+\mu(W-I_h(W))-\mu W.
\]
Take $L^2-$inner product with $W$ on both sides of the equation, we have
\begin{equation}\label{eq:prop1 ineq1}
\begin{split}
 \frac{1}{2}\frac{d}{dt}||W||^2_2 &= \lan G-\tG,W \ran+\lan S(U)-S(V),W \ran + \mu\lan W-I_h(W),W \ran-\mu||W||_2^2\\
 &\leq ||G-\tG||_2||W||_2+||S(U)-S(V)||_2||W||_2+\mu||W-I_h(W)||_2||W||_2-\mu||W||_2\\
 &\leq ||G-\tG||_2||W||_2+L_s||W||^2_2+\mu||W-I_h(W)||_2||W||_2-\mu||W||^2_2\\
 & = ||G-\tG||_2||W||_2+\mu||W-I_h(W)||_2||W||_2-(\mu-L_s)||W||^2_2.
\end{split}
\end{equation}
To bound the error $||G-\tG||_2$, we employ the condition \eqref{eq:observ cond2} and write
\begin{equation}\label{eq:prop1 ineq2}
\begin{split}
G-\tG &= I_h(G)-\tG\\
& = I_h\big(\partial_t U+\partial_x U-S(U)\big)-\big(\partial_tI_h(U)+I_h(\partial_x\tU-S(\tU)\big)\\
& = I_h(\partial_x U-\partial_x \tU)-I_h(S(U)-S(\tU)).
\end{split}
\end{equation}
The temporal derivatives are canceled due to $I_h \partial_t(\cdot) = \partial_t I_h(\cdot)$. Subsequently, using conditions \eqref{eq:Ih cond 2}, \eqref{eq:Ih cond 3} leads to
\[
\begin{split}
||G-\tG||_2& \leq ||I_h\big(\partial_x (U-\tU)\big)||_2+||I_h\big(S(U)-S(\tU)\big)||_2\\
&\leq (c_1 h^{-1}+C_0 L_s)||U-\tU||_2\\
& = (c_1 h^{-1}+C_0 L_s)||W-I_h(W)||_2.
\end{split} 
\]
By substituting the above inequality into \eqref{eq:prop1 ineq1} and rearranging, we obtain 
\begin{equation}\label{eq:prop1 ineq3}
\setlength{\jot}{10pt}
\begin{split}
\frac{d}{dt}||W||_2+(\mu-L_s)||W||_2 \leq (c_1h^{-1}+C_0 L_s+\mu)||W-I_h(W)||_2.
\end{split}
\end{equation}
Applying Gronwall's inequality and employing condition \eqref{eq:observ cond1} yields
\[
\setlength{\jot}{10pt}
\begin{split}
||W(t)||_2&\leq \big\{||W(0)||_2+(c_1h^{-1}+C_0 L_s+\mu)\int^t_0 ||W(\tau)-I_h\big(W(\tau)\big)||_2 \, e^{(\mu-L_s)\tau}d\tau\}  \, e^{-(\mu-L_s)t}\\
&\leq ||W(0)||_2 \, e^{-(\mu-L_s)t}+\frac{(c_1h^{-1}+C_0 L_s+\mu)C_h}{\mu-L_s-\eps} \, e^{-\eps t}.
\end{split}
\]
The desired estimate \eqref{eq:state exp decay} is recovered. \hfill $\square$
\medskip

\begin{remark}
To get a better understanding of the convergence of the algorithm, especially to understand the condition \eqref{eq:observ cond1}, it is helpful to look at a simpler example:
\begin{subequations}\label{eq: reduced ode}
\begin{equation}
\partial_tU(x,t) = G(x,t)+s U(x,t),
\end{equation}
and the corresponding nudged system becomes:
\begin{equation}
\left\{\begin{array}{l}
    \partial_t V = \tG+ sV+\mu(I_h(U)-I_h(V)),  \\
    \\
    \tG(t) = I_h(\partial_t U)-sI_h(\tU),\quad \tU = I_h(U)+(I-I_h)(V),
\end{array}
\right.
\end{equation}
\end{subequations}
where $s\in\R$ is a constant. We consider the setting where $I_h$ is a compact symmetric operator on $L^2(\R^d)$. Then there is an orthornormal basis of $L^2(\R^d)$ consisting of eigenfunctions, $p_1, p_2, \cdots$, of $I_h$, with corresponding eigenvalues $\lambda_1,\lambda_2,\cdots$. Let us assume $c_0 = 1$ in \eqref{eq:Ih cond 2} and $I_h$ is semi-definite, so that $0\leq \lambda_i\leq 1$. 

Now, we can write all quantities using an eigenfunction basis as
\[
U(x,t) = \sum_i u_i(t) p_i(x), \ \ G(t) = \sum_i g_i(t) p_i(x),\ \
V(t) = \sum_i v_i(t) p_i(x), \ \ \tG(t) = \sum_i \tg_i(t)p_i(x).
\]
Projection of equations \eqref{eq: reduced ode} onto the eigenfunction modes yields
\begin{equation*}
\left\{\begin{array}{l}
u_i'(t) = g_i(t) + s u_i(t),\\
\\
v_i'(t) = \tilde g_i(t) + s v_i(t) + \mu \lambda_i (u_i(t) - v_i(t)),\\
\\
\tilde u_i(t) = \lambda_i u_i(t) + (1-\lambda_i)v_i(t),\\
\\
\tilde g_i(t) = \lambda_i(g_i + s u_i(t)) - s\lambda_i^2 u_i(t) - s\lambda_i(1-\lambda_i)v_i(t).
\end{array}\right.
\end{equation*}
Taking subtraction $u'_i(t)-v'_i(t)$ and using that $\lambda_i g_i = g_i$ (which follows from \eqref{eq:observ cond2}) yields
\[\frac{d}{dt} (v_i(t) - u_i(t)) = \left( s - \mu \lambda_i - s\lambda_i(1-\lambda_i)\right)(v_i(t) - u_i(t)).\]
Therefore,
\[v_i(t) - u_i(t) = e^{ (s - \mu \lambda_i - s\lambda_i(1-\lambda_i))t}(v_i(0) - u_i(0)),\]
and so 
\begin{equation}\label{eq:ode converge}
V(t) = U(t) + \sum_{i=1}^{\infty}e^{ (s - \mu \lambda_i - s\lambda_i(1-\lambda_i))t}(v_i(0) - u_i(0))p_i.
\end{equation}
This implies an interplay between the structure of source terms and the conditions for convergence.

If $s>0$, to guarantee exponential convergence for any initial condition we cannot have any $\lambda_i=0$, or if we do, then we must ask that $u_i(0)=v_i(0)$. This is reasonable to expect -- we cannot damp any of the error modes that belong to the kernel of $I_h$. The lower bound of $\mu$ becomes
\begin{equation}\label{eq:mu_cond}
\mu\geq \frac{s+\eps}{\lambda_i}, \quad \eps>0,
\end{equation}
for all $i$ such that $u_i(0)\neq v_i(0)$, from which we obtain $e^{-\eps t}$ convergence. This is connected to condition \eqref{eq:observ cond1}, since we can compute
\[
||W(t)-I_h\big(W(t)\big)||^2_2 = \sum_{i=1}^{\infty} (\lambda_i - 1)^2 (v_i(0) - u_i(0))^2e^{ 2(s - \mu \lambda_i - s\lambda_i(1-\lambda_i))t}.
\]
Suppose that \eqref{eq:mu_cond} holds whenever $0 < \lambda_i < 1$ and $u_i(0) \neq v_i(0)$, and $u_i(0)=v_i(0)$ whenever $\lambda_i=0$. Then we recover  
\[
||W(t)-I_h\big(W(t)\big)||_2\leq ||W(0)||_2 \, e^{-\eps t}.
\]
In other words, condition \eqref{eq:observ cond1} may be regarded as an implication from an appropriate initial state and a sufficiently large $\mu$.

The situation becomes quite different when $s<0$.
In this scenario, the source term has stabilizing effect on the dynamics and helps dissipate the error. Note that in \eqref{eq:ode converge} the exponents are bounded by
\[
s-\mu\lambda_i-s\lambda_i(1-\lambda_i)\leq s-\lambda_i(1-\lambda_i)s \leq \frac{3}{4}s<0.
\]
This implies unconditional convergence regardless of $\mu$, and \eqref{eq:observ cond1} is automatically recovered with $\displaystyle \eps = -\frac{3}{4}s > 0$. However, if the spatial transport is also considered, the dissipation of the source term is not necessarily sufficient to absorb the error from the convection term. Observational feedback control is still needed to enhance the convergence in general.
\end{remark}

The convergence of state further implies the convergence of force, as stated in the following proposition. 

\begin{proposition}\label{prop:force recovery}\textbf{\emph{(Convergence of the force)}}
Let $(U,G)$ and $(V,\tG)$ be the strong solutions of the hydrodynamics \eqref{eq:hydro} and the nudged system \eqref{eq:nudge} associated with $I_h$ satisfying conditions \eqref{eq:Ih conditions}. Assume that the convection flux $F(\cdot)$ and the source term $S(\cdot,\bx,t)$ are Lipschitz continuous with Lipschitz constants $L_f$, $L_s>0$. If the nudged solution $V(\bx,t)$ recovers the true data $U(\bx,t)$ asymptotically in time, $||U(t)-V(t)||_2 \rightarrow 0$ as $t\rightarrow\infty$, then
\[
||I_h(G(t))-\tG(t)||_2 \rightarrow 0, \quad t\rightarrow\infty.
\]
In particular, if the force function is completely observable, $I_h(G) = G$, then 
\[
||G(t)-\tG(t)||_2 \rightarrow 0, \quad t\rightarrow\infty.
\]
\end{proposition}
\emph{Proof}. We have
\[
\tG = \partial_t I_h(U) + I_h(\nabla \cdot F(\tU)-S(\tU)), \quad \tU = I_h(U)+(I-I_h)(V),
\]
and
\[
I_h(G) = \partial_t I_h(U) + I_h(\nabla \cdot F(U)-S(U)).
\]
The difference between the above equations yields
\[
\begin{split}
    I_h(G)-\tG = I_h(\nabla\cdot(F(U)-F(\tU)))-I_h(S(U)-S(\tU)).
\end{split}
\]
Using the estimates \eqref{eq:Ih cond 2}, \eqref{eq:Ih cond 3}, it follows that
\begin{equation}\label{eq:prop2 ineq1}
\begin{split}
    ||I_h(G)-\tG||_2 &\leq ||I_h(\nabla\cdot(F(U)-F(\tU)))||_2+||I_h(S(U)-S(\tU))||_2\\
    &\leq c_1h^{-1} ||F(U)-F(\tU)||_2+C_1 h ||S(U)-S(\tU)||_2\\
    &\leq (c_1 h^{-1}L_f+C_1 h L_s)||U-\tU||_2  \\
    & = (c_1 h^{-1}L_f+C_1 h L_s)||(U-V)-I_h(U-V)||_2 \\
    & \leq (c_1 h^{-1}L_f+C_1 h L_s)C_0||U-V||_2.
\end{split}
\end{equation}
Hence, the convergence of the state, $||U(t)-V(t)||_2 \rightarrow 0$, implies the convergence of the reconstructed force, $||I_h(G(t))-\tG(t)||_2\rightarrow 0$. 
\hfill$\square$

\medskip
In contrast to the earlier discussions on the incompressible Navier-Stokes equations \cite{azouani2014continuous,gesho2016computational}, our analysis in Proposition \ref{prop:state recovery} employs the extra condition \eqref{eq:observ cond1}. Essentially, the purpose is to ensure that the relaxation term $\mu (I_h(U)-I_h(V))\overset{t\rightarrow\infty}{\longrightarrow} \mu(U-V)$. The difficulty in weakening this condition arises from the lack of viscosity, which, if exists, would help dissipate the interpolation error. For instance, by adding viscous terms  $\nu\partial^2_x U$ and $\nu\partial^2_x V$, with $\nu>0$, to equations \eqref{eq:linear eqn} and \eqref{eq:linear nudge}, we can derive the following inequality along the lines of \eqref{eq:prop1 ineq1} and \eqref{eq:prop1 ineq2},
\[
\setlength{\jot}{10pt}
\begin{split}
\frac{1}{2}\frac{d}{dt}||W||^2_2+(\mu-L_s)||W||^2_2 &\lesssim \mu ||W-I_h(W)||_2||W||_2-\nu ||\partial_x W||^2_2\\
& \lesssim c_1\mu h ||w||_{H^1}||w||_2-\nu||\partial_x w||^2_2\\
& \lesssim \frac{\mu}{2}||w||^2_2 +\frac{\mu c_1^2 h^2}{2}||\partial_x w||^2_2-\nu||\partial_x w||^2_2.
\end{split}
\]
The interpolation error is absorbed by the viscosity for $L_s \lesssim \mu\lesssim \nu/h^2$, thereby ensuring the exponential convergence. This is similar to the results in \cite{azouani2014continuous,gesho2016computational}.

Back to the setting of inviscid hydrodynamics, \eqref{eq:Ih cond 3} implies $||w-I_h(w)||_2 \sim h||\partial_x w||_2$, hence the condition \eqref{eq:observ cond1} amounts to $\displaystyle ||\partial_x w(t)||_2\lesssim \frac{1}{h} e^{-\eps t}$. This indicates that when the data has strong oscillations (large $\dot{H}^1-$norm), the analytical convergence may not be guaranteed under insufficient observations. Here we prove the worst-case scenario with respect to generic bounded data.

\begin{proposition}\label{prop:error for sparse observ}\textbf{\emph{(Error estimate under insufficient observations)}}
Let $(U,G)$ be a strong solution of the hydrodynamics \eqref{eq:hydro} driven by the Lipschitz continuous convection flux and source term. Let $(V,\tG)$ be a strong solution of the nudged system \eqref{eq:nudge}, where the associated $I_h$ satisfies conditions \eqref{eq:Ih conditions}. Assume that the data are uniformly bounded in $H^1(\R^d)$, i.e., 
\begin{equation*}
||G(t)||_{H^1}, \hspace{0.3em} ||\tG(t)||_{H^1} \leq M_g, \quad ||U(t)||_{H^1}, \hspace{0.3em} ||V(t)||_{H^1} \leq M_u,
\end{equation*}
for all $t>0$, with $M_u$, $M_g>0$. Then the nudged solutions approach the true solution time-asymptotically:
\begin{subequations}\label{eq:worst est}
\begin{equation}\label{eq:worst est state}
\overline{\lim_{t\rightarrow\infty}}||U(t)-V(t)||_2 \leq \frac{\kappa_1}{\mu-L_s}, 
\end{equation}
\begin{equation}\label{eq:worst est force}
\overline{\lim_{t\rightarrow\infty}}||G(t)-\tG(t)||_2 \leq \frac{\kappa_2}{\mu-L_s}+C_1M_g h,
\end{equation}
\end{subequations}
where $\kappa_1 = 2(L_f M_u+M_g+C_1M_uh)$, $\kappa_2 = (c_1 h^{-1}L_f+C_1hL_s)C_0\kappa_1$.

\end{proposition}
\emph{Proof.} Denote $W = U-V$, the difference of \eqref{eq:hydro} and \eqref{eq:nudge eqn} yields
\[
\frac{\partial W}{\partial t}+\nabla\cdot(F(U)-F(V)) = G-\tG+S(U)-S(V)+\mu(W-I_h(W))-\mu W.
\]
Taking $L^2-$inner product with $W$ leads to
\[
\setlength{\jot}{10pt}
\begin{split}
 \frac{1}{2}\frac{d}{dt}||W||^2_2 &=  -\lan W, \nabla\cdot\big(F(U)-F(V)\big)\ran+\lan W,G-\tG\ran+\lan W,S(U)-S(V)\ran+\lan W,W-I_h(W)\ran-\mu||W||^2_2\\
 & \leq \int_{\R^d} (F(U)-F(V)):\nabla W dx +||G-\tG||_2||W||_2+L_s||W||^2_2+||W||_2||W-I_h(W)||_2-\mu||W||^2_2\\
 &\leq L_f ||W||_2||W||_{H^1}+2M_g||W||_2+L_s||W||^2_2+C_1 h ||W||_2||W||_{H^1}-\mu||W||^2_2\\
 &\leq 2(L_f M_u+M_g+C_1  M_u h)||W||_2-(\mu-L_s)||W||^2_2.
\end{split}
\]
Rearrange the above inequality, we obtain
\[
\frac{d}{dt} ||W||_2+(\mu-L_s)||W||_2\leq 2(L_f M_u+M_g+C_1  M_u h).
\]
Integration yields
\[
||W(t)||_2 \leq ||W(0)||_2 e^{-(\mu-L_s)t}+\frac{2(L_f M_u+M_g+C_1M_uh)}{\mu-L_s}(1-e^{-(\mu-L_s)t}),
\]
the bound \eqref{eq:worst est state} is recovered by passing the limit $t\rightarrow\infty$. Furthermore, 
\[
\setlength{\jot}{10pt}
\begin{split}
||G(t)-\tG(t)||_2 &\leq ||I_h(G(t))-\tG(t)||_2+||G(t)-I_h(G(t))||_2\\
& \leq (c_1 h^{-1}L_f+C_1hL_s)C_0||U(t)-V(t)||_2+C_1h||G(t)||_{H^1}\\
& \leq (c_1 h^{-1}L_f+C_1hL_s)C_0||W(t)||_2+C_1 M_g h.
\end{split}
\]
The second line is obtained using the inequality \eqref{eq:prop2 ineq1}. The error bound of force recovery \eqref{eq:worst est force} is a direct consequence of the large-time limit.\hfill $\square$

\begin{remark}
It is noticed that the error bounds in Propositions \ref{prop:state recovery}, \ref{prop:force recovery}, and \ref{prop:error for sparse observ} involve $h^{-1}$ dependence. Indeed, along the lines of \eqref{eq:prop2 ineq1}, the $h^{-1}$ term occurs when estimating the error of force recovery, $||\tG-I_h(G)||_2$, where the usage of interpolation bound \eqref{eq:Ih cond 2} yields 
\[
||I_h\big(\nabla\cdot(F(U)-F(\tU))\big)||_2\leq c_1h^{-1}||F(U)-F(\tU)||_2.
\]
Nevertheless, assuming the data admits sufficiently weak oscillations such that
\[
\sup_{t\geq 0} \frac{||\nabla\cdot(F(U)-F(\tU))||_2}{||F(U)-F(\tU)||_2}\leq C<\infty,
\]
and the above estimate can be modified to
\[
||I_h\big(\nabla\cdot(F(U)-F(\tU))\big)||_2\leq c_0C||F(U)-F(\tU)||_2,
\]
following the lines of \eqref{eq: Ih weak oscillation bound}. Subsequently, the $h^{-1}$ dependence in the error bounds will be removed.

Despite the difficulty in deriving an optimal error bound under less ideal assumptions, in actual computations, we should still expect a more accurate solution under a finer observation. 
Meanwhile, though the error bound in Proposition \ref{prop:error for sparse observ} does not offer an explicit criterion for determining the minimum observation density required to meet a specific error threshold, we will demonstrate that our proposed algorithm maintains reasonable accuracy even with relatively sparse data through numerical experiments in Section \ref{subsec:noisy}.
\end{remark}

\subsection{Discretization of nudged system}
To solve the nudged system \eqref{eq:nudge}, since both the state $V$ and the force $\tG$ are evolving in time, we introduce the following discretization. Given the solutions $V^n = \{V^n_j\}^N_{j=1}$ and $\tG^n = \{\tG^n_j\}^N_{j=1}$ at time level $t = t^n$, the solution at $x = x_j$ and $t = t^{n+1}$ are computed through:
\[
\left\{\begin{array}{l}\displaystyle
\frac{V_j^{n+1}-V_j^n}{\Delta t}+\delta^{\Dx}_x F(V^n)_j = \tG^n_j+S(V^n,x_j,t^n) \,,\\
\\ \displaystyle
 \tG^{n+1}_j = \partial_t I_h(U(t^{n+1}))_j+I_h(\delta^{\Dx}_x F(\tU^{n+1})-S(\tU^{n+1}))_j\,,\\
\\ \displaystyle
\tU^{n+1} = I_h(U(t^{n+1}))+(I-I_h)(V^{n+1})\,.
\end{array}\right.
\]
Here $\Dt$ is the time step, and $\delta^{\Dx}_x F$ is the approximate flux derivative obtained with non-oscillatory discretization such as ENO/WENO schemes with spatial mesh size $\Dx$. We assume that the observation $I_h(U(t))$ is recorded continuously in time, whereas $\partial_t I_h(U(t^{n+1}))$ can be computed with high order finite difference methods. As such, the above forward Euler discretization can be directly replaced by higher order strong stability preserving (SSP) Runge-Kutta methods \cite{gottlieb2001strong}. 

It is worth noting that the numerical scheme described above yields nodal values, $\{V_j\}$, over the computational grids, $\{x_j\}$. The interpolation $I_h(\cdot)$, however, is constructed based on the data over the observation grids, $\{x^{ob}_j\}$, which may not coincide with ${x_j}$ in general. Therefore, appropriate interpolation techniques (such as linear interpolation or cubic spline) may be necessary to approximate nodal values ${(x^{ob}_j, V^{ob}_j)}$. These approximations are then used to compute the interpolations ${(x_j, I_h(V)_j)}$.

We also emphasize the simplicity of our discretization. In the previous works \cite{martinez2022reconstruction,farhat2024identifying}, it was proposed to generate a \emph{sequence} of approximations, $(V_k(\bx,t), G_k(\bx,t))|_{t\geq t_k}$, $t_k\gg t_{k-1}$, $k = 0,1,2,\cdots$. The solution $V_k(\bx,t)$ was obtained by solving the nudged equations with the re-initialization, $V_k(\bx,t_{k}) = V_{k-1}(\bx,t_{k})$, coupled with the frozen force, $G_{k-1}(\bx,t)|_{t\geq t_{k}}$, over the time interval $t\in[t_k,\infty)$. Although this algorithm guarantees exponential convergence when applied to the incompressible Navier-Stokes equations, the need to repeatedly solve the equation can be infeasible in industrial applications. In contrast, our discretization method updates the states and the force simultaneously \emph{on the fly}, hence significantly reducing the computational cost.


\subsection{Smooth fitting with kernel regression}
\label{sec:kernel regression}
In practical computations, the accuracy of approximate solutions can depend on the selection of $I_h$, particularly in scenarios with restricted observational data. Common methods for achieving a smooth interpolation of observed data involve using Lagrangian interpolation with nodal values distributed uniformly at a mesh size $h$, or projecting onto truncated Fourier series with a cut-off frequency of $1/h$. However, these approaches using high-order polynomials or trigonometric polynomials may encounter significant challenges such as numerical oscillations at boundaries (known as Runge's phenomenon) or sharp extrema (referred to as Gibbs' phenomenon).   

To reduce the undesired numerical artifacts, we propose to \emph{fit} the observed data with kernel regression,
\begin{equation}\label{eq:KDE}
I_h(U)(\bx) = \frac{\sum^{N_{ob}}_{j=1}K_{\sig}(|\bx-\bx^{ob}_j|)U^{ob}_{j}}{\sum^{N_{ob}}_{j=1}K_{\sig}(|\bx-\bx^{ob}_j|)},
\end{equation}
where $U^{ob}_j$ is the observation at $\bx^{ob}_j\in\R^d$, $K_\sig(\cdot)$ is the Gaussian kernel with bandwidth $\sig$,
\begin{equation}\label{eq:gaussian kernel}
K_\sig(|\bx|) = \frac{1}{\sig^d}\exp\{-\frac{|\bx|^2}{2\sig^2}\}.
\end{equation}
A larger $\sig$ leads to stronger kernel smoothing, whereas a smaller $\sig$ generates a theoretically more accurate fitting. However, taking $\sig$ too small in actual computations can result in severe arithmetic errors. Hence, to achieve a balance between solution sharpness and computational stability, we suggest taking $\displaystyle\frac{\sqrt{d}}{2}h<\sig\leq h$. It's important to note that the approximate data  $I_h(U)$ is derived from weighted averages of nearby observations, which generally do not preserve nodal values at interpolation points. In practical applications where observations may be noisy due to measurement errors, using kernel regression instead of classical interpolation methods can better mitigate over-fitting issues by appropriately adjusting the kernel bandwidth.

A common challenge when using kernel regression \eqref{eq:KDE} is the larger error observed at the boundary of the compact support. This is due to the boundary conditions of the data $ U(\bx)$ not being incorporated, leading to higher accuracy in the interior of the domain compared to near the boundary. To address this, we impose a simple boundary correction via ghost cell extensions. For the convenience of demonstration, we consider uniformly distributed observation grids $\{x^{ob}_j\}^{N_{ob}}_{j=1}$ over the computational interval $[a,b]$ with $x^{ob}_1 = a$ and $x^{ob}_{N_{ob}} = b$. We modify the kernel regression  with
\begin{equation}\label{eq:KDE modified}
 I_h(U)(x) = \frac{\sum^{N_{ob}+m}_{j=-(m-1)}K_{\sig}(|x-x^{ob}_j|)U^{ob}_{j}}{\sum^{N_{ob}+m}_{j=-(m-1)}K_{\sig}(|x-x^{ob}_j|)}.
\end{equation}
The ghost nodes are given by
\[
x^{ob}_{1-j} = a-jh, \quad x^{ob}_{N_{ob}+j} = b+jh, \quad j = 1,2,\cdots,m.
\]
This way, boundary conditions are encoded implicitly in the extrapolation of
$\{U^{ob}_{j}\}_{j\leq 0}$ and $\{U^{ob}_j\}_{j\geq N_{ob}+1}$. For instance, the extended data at the left boundary are set to be
\[
U^{ob}_{1-j} = \left\{\begin{array}{ll}
  U^{ob}_{N_{ob}-j} &\text{, periodic boundary}  \\
  &\\
  2U_a-U^{ob}_{1+j}  &\text{, Dirichlet boundary $U(a) = u_a$}\\
  &\\
  U^{ob}_{1+j}  &\text{, Neumann boundary $U'(a) = 0$}
\end{array}\right., \quad j = 1,2,\cdots, m.
\]
The data can be extrapolated at the right boundary in a similar manner. 
Although the kernel regression does not meet the conditions specified in \eqref{eq:Ih conditions}, it does satisfy estimates similar to those in \eqref{eq:Ih cond 3}, as is verified in the following result.
\begin{theorem}\label{thm:kernel_regression}
Let $U^{ob}_{j} = U(\bx^{ob}_j)$, where $U:\Omega\to \R$, and assume that $K_\sig(\bx)=0$ for $|\bx|\geq C\sig$. Then the following hold. 
\begin{enumerate}[(i)]
\item If $U$ is Lipschitz continuous with Lipshcitz constant $L$, then 
\[\|I_h(U) - U\|_\infty \leq CL\sig.\]
\item If $\nabla U$ is Lipschitz continuous with Lipschitz constant $L$, then 
\[\|I_h(U) - U\|_2 \leq  C \sig \|U\|_{H^1} + LC^2\sigma^2|\Omega|^{1/2}.\]
\end{enumerate}
\end{theorem}
\emph{Proof.} For simplicity, let us write 
\[d_\sig(x) = \sum^{N_{ob}}_{j=1}K_{\sig}(|\bx-\bx^{ob}_j|).\]
Then we have 
\begin{align*}
|I_h(U)(\bx) - U(\bx)|  &= \left|\frac{1}{d_\sigma(x)}\sum^{N_{ob}}_{j=1}K_{\sig}(|\bx-\bx^{ob}_j|)U(\bx^{ob}_{j}) - U(\bx) \right|\\
&\leq  \frac{1}{d_\sig(x)}\sum^{N_{ob}}_{j=1}K_{\sig}(|\bx-\bx^{ob}_j|)\left|U(\bx^{ob}_{j}) - U(\bx)\right| \\
&\leq  \frac{L}{d_\sig(x)}\sum^{N_{ob}}_{j=1}K_{\sig}(|\bx-\bx^{ob}_j|)\left|\bx^{ob}_{j} - \bx\right|\\
&\leq \frac{CL\sig}{d_\sig(x)}\sum^{N_{ob}}_{j=1}K_{\sig}(|\bx-\bx^{ob}_j|) = CL\sig,
\end{align*}
which establishes (i).

For (ii), we continue the computation above to obtain
\begin{align*}
|I_h(U)(\bx) - U(\bx)|  &\leq  \frac{1}{d_\sig(x)}\sum^{N_{ob}}_{j=1}K_{\sig}(|\bx-\bx^{ob}_j|)\left|U(\bx^{ob}_{j}) - U(\bx)\right| \\
&\leq  \frac{1}{d_\sig(x)}\sum^{N_{ob}}_{j=1}K_{\sig}(|\bx-\bx^{ob}_j|)\left(\left| \nabla U(\bx)\cdot(\bx^{ob}_j - \bx)\right| + L|\bx^{ob}_j - \bx|^2\right) \\
&\leq  \frac{1}{d_\sig(x)}\sum^{N_{ob}}_{j=1}K_{\sig}(|\bx-\bx^{ob}_j|)\left(\left| \nabla U(\bx)\right|C\sig + LC^2\sig^2\right)\\
&\leq \left| \nabla U(\bx)\right|C\sig + LC^2\sig^2.
\end{align*}
The proof is completed by taking the $L^2$ norm on both sides. \hfill $\square$

\begin{remark}\label{rem:der_est}
The estimates in Theorem \ref{thm:kernel_regression} are similar to \eqref{eq:Ih cond 3} with $m=1$, provided we use that $\sigma \sim h$, though they are not exactly the same. Unfortunately, it is not possible to prove estimates similar to \eqref{eq:Ih cond 2}; instead, we can prove estimates on the derivatives of the interpolated function, that is the quantity $\|\partial^\alpha I_h(U)\|$ is straightforward to estimate, since in the case that $K_\sig(\bx) = \frac{1}{\sig^d}K_1(\bx/\sig)$ we have
\begin{equation}\label{eq:KDE_diff}
\partial^\alpha I_h(U)(\bx) = \frac{\sum^{N_{ob}}_{j=1}\sig^{-|\alpha|}\partial^\alpha K_{\sig}(|\bx-\bx^{ob}_j|)U^{ob}_{j}}{\sum^{N_{ob}}_{j=1}K_{\sig}(|\bx-\bx^{ob}_j|)}.
\end{equation}
However, such estimates are not helpful in our analysis, so we do not pursue them further. Even though kernel regression does not fully satisfy \eqref{eq:Ih conditions}, the method still produces satisfactory numerical solutions, as demonstrated by the examples in the next section. 
\end{remark}

In summary, our construction of interpolant operator via kernel regression provides significant practical benefits in the following aspects:

\noindent(i) Convenient implementation of flexible boundary conditions: kernel regression allows for straightforward handling of complex and non-standard boundary conditions, which is often a crucial challenge suffered by other interpolants, e.g. spectral projections or piecewise linear interpolation;\\
\noindent(ii) Robust processing of noisy observations: our approach demonstrates inherent resilience to noisy observational data, facilitating more reliable interpolation in real-world scenarios where measurements are often imperfect. This is explicitly shown in the tests with noisy data presented in the following section \ref{subsec:noisy}.

These advantages are critical for potential extensions to industrial applications, such as reactor simulations with complex geometry and possibly inaccurate measured data. In simpler computational settings, alternative choices of the interpolant operator may also provide satisfactory results. As the main focus of our paper is to validate the nudging-based feedback control framework as applied to generic inviscid hydrodynamics, particularly when recovering moments from incomplete observations (see Section \ref{sec:moment-recovery}), further investigations of other interpolants will be left for future study.


\section{Numerical experiments I -- states and force recovery}
\label{sec:numerical1}

We verify the effectiveness of the algorithm \eqref{eq:nudge} in both scalar and system test cases. In all experiments, unless stated otherwise, the operator $I_h$ is derived using the kernel regression \eqref{eq:KDE modified} with $m = 3$ ghost nodes applied at each boundary. The kernel bandwidth in equation  \eqref{eq:gaussian kernel} is set to $\sig = h$. The observation grids $\{x^{ob}_j\}$ are uniformly placed over the computational domain. The nudged equations are discretized with the 5th-order WENO finite difference scheme in space, and time integration is performed with an explicit SSPRK3 method with $CFL = 0.7$. Reference solutions are obtained by discretizing the exact dynamics with WENO5-SSPRK3 on fine meshes.

\subsection{Scalar test case}\label{sec:scalar test}
Consider the following exact dynamics: 
\begin{equation}\label{eq:scalar dynamics}
\partial_t u + \partial_x u = -\partial_x p + S(u) , \quad x\in[0, 2\pi],
\end{equation}
with 
\[
p = \frac{u^2}{6}, \quad S(u) = 0.2\sqrt{1+u^2}.
\]
The initial data consists of sine and cosine waves with different frequencies,
\[
u_0(x) = -0.8\sin(x)+0.4\sin(2x)+0.02\cos(10x).
\]
Periodic boundary conditions are applied at both endpoints. In this test case, the pressure gradient $g = -\partial_x p$ generates high-frequencies and sharp extrema over time due to the non-linearity.

Assume that the source function $S(\cdot)$ is known, our goal is to recover the  force term, $\displaystyle g = -\partial_x (\frac{u^2}{6})$, by evolving the nudged system
\begin{equation}\label{eq:scalar nudging}
\left\{\begin{array}{l}
\partial_t v +\partial_x v = \tg+0.2\sqrt{1+v^2}+\mu(I_h(u)-I_h(v))\,, \\
\\
\tg = \partial_t I_h(u)+I_h(\partial_x\tu-0.2\sqrt{1+\tu^2}), \quad \tu = I_h(u)+(I-I_h)(v)
\end{array}\,,
\right.
\end{equation}
with zero initial data,
\[
\begin{split}
v_0(x) = 0, \quad \tg_0(x) = 0.
\end{split}
\]
The nudging coefficient $\mu$ is set to be $3$. We intend to recover the solution state $u$ and the force function $g$ from $N_{ob} = 150$ spatial observation points. The nudged system is discretized with $N = 800$ equidistant grid points.

\begin{figure}[h!]
    \centering
    \begin{subfigure}{0.3\textwidth}
    \includegraphics[scale = 0.22]{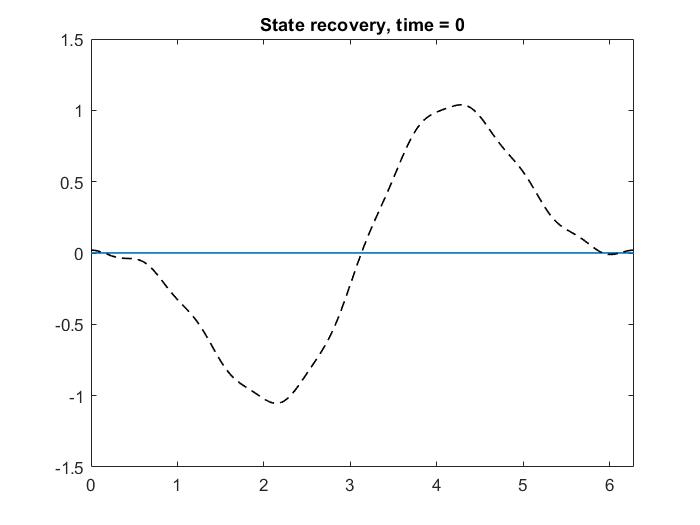}
    \subcaption{$t = 0$}
    \end{subfigure}
    \quad
    \begin{subfigure}{0.3\textwidth}
    \includegraphics[scale = 0.22]{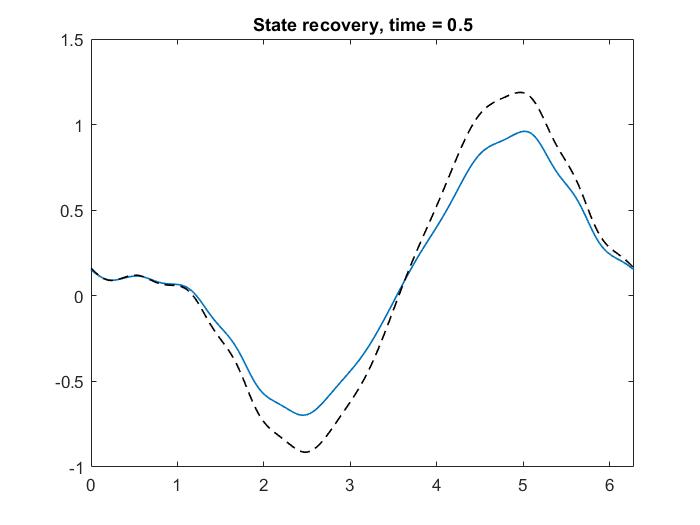}
    \subcaption{$t = 0.5$}
    \end{subfigure}
    \\
    \begin{subfigure}{0.3\textwidth}
    \includegraphics[scale = 0.22]{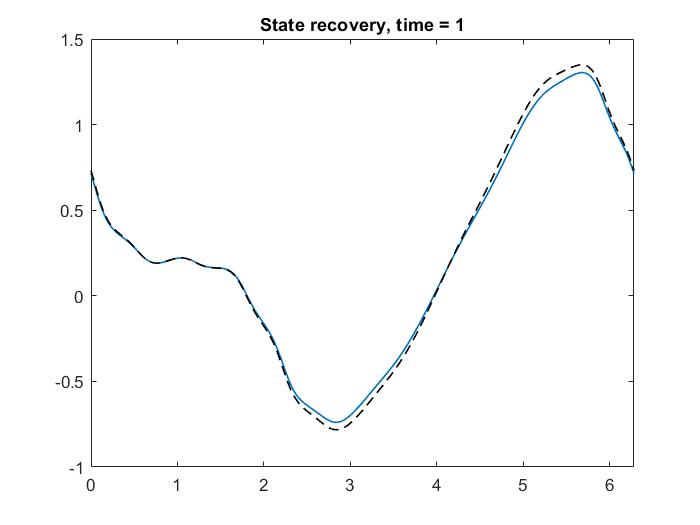}
    \subcaption{$t = 1$}
    \end{subfigure}
    \quad
    \begin{subfigure}{0.3\textwidth}
    \includegraphics[scale = 0.22]{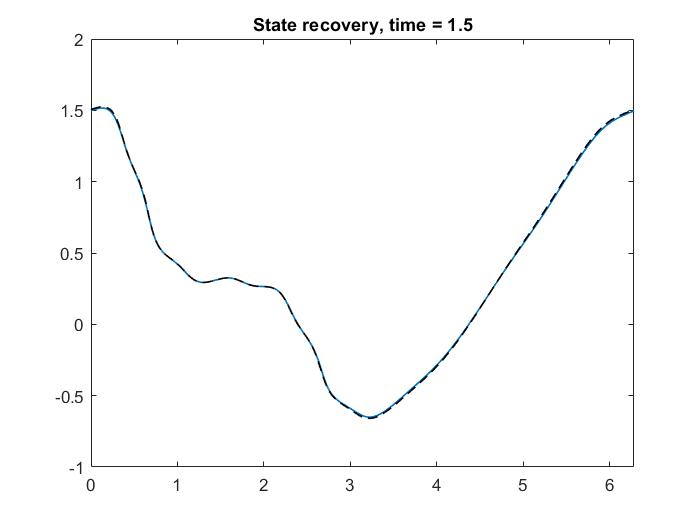}
    \subcaption{$t = 1.5$}
    \end{subfigure}
    \caption{Scalar dynamics. State recovery. Nudged solution (blue lines) vs. reference solution (black dash lines).} 
    \label{fig:linear state}
\end{figure}

\begin{figure}[h!]
    \centering
    \begin{subfigure}{0.3\textwidth}
    \includegraphics[scale = 0.22]{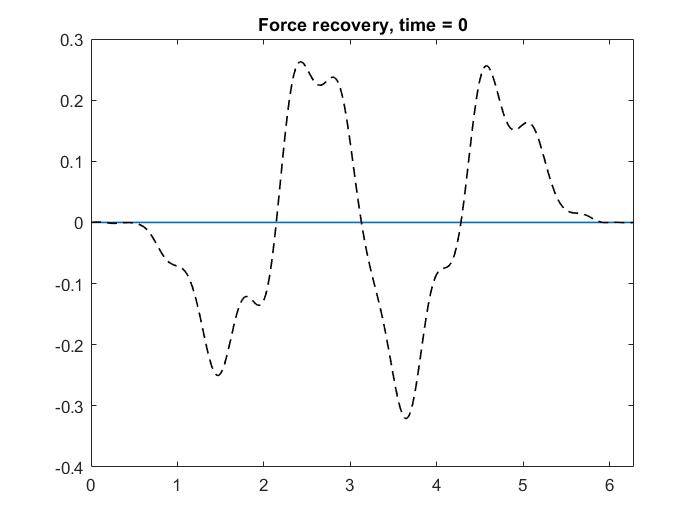}
    \subcaption{$t = 0$}
    \end{subfigure}
    \quad
    \begin{subfigure}{0.3\textwidth}
    \includegraphics[scale = 0.22]{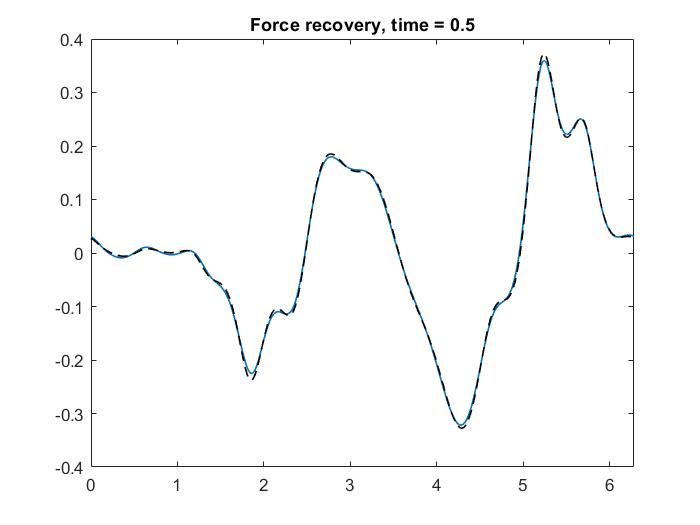}
    \subcaption{$t = 0.5$}
    \end{subfigure}
    \\
    \begin{subfigure}{0.3\textwidth}
    \includegraphics[scale = 0.22]{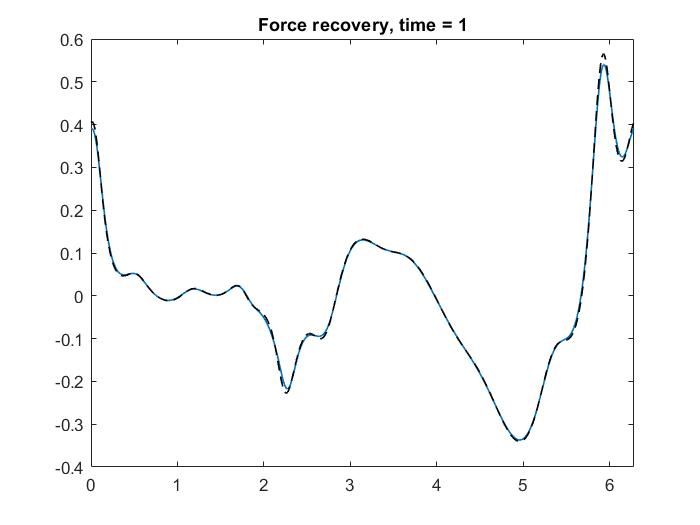}
    \subcaption{$t = 1$}
    \end{subfigure}
    \quad
    \begin{subfigure}{0.3\textwidth}
    \includegraphics[scale = 0.22]{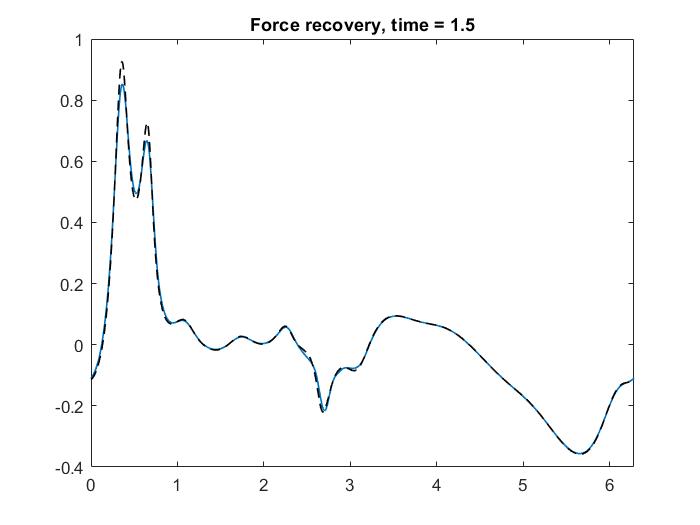}
    \subcaption{$t = 1.5$}
    \end{subfigure}
    \caption{Scalar dynamics. Force recovery. Nudged solution (blue lines) vs. reference data (black dash lines).} 
    \label{fig:linear force}
\end{figure}

\begin{figure}[h!]
    \centering
    \includegraphics[scale = 0.23]{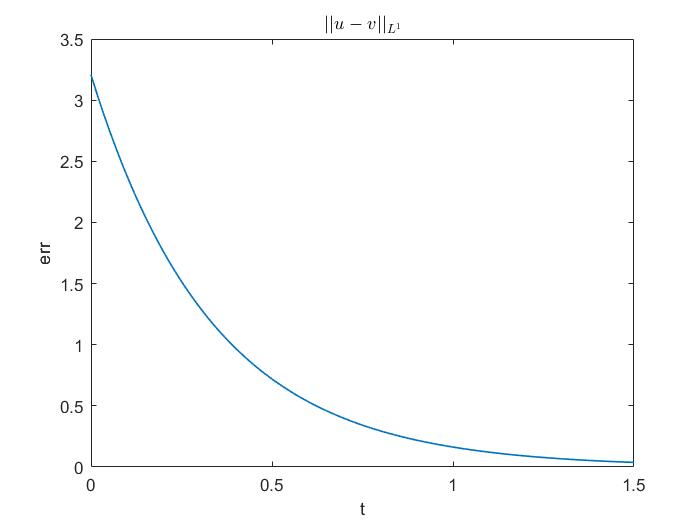}
    \caption{Scalar dynamics. History of $L^1$ state error. }
    \label{fig:linear err}
\end{figure}

Figures \ref{fig:linear state}, \ref{fig:linear force} display the evolution of approximate solutions up to $T = 1.5$. The computed solutions $v$ and $\tilde g$ converge rapidly to the reference solutions $u$ and $g$ over time. Although the smoothing effect of kernel regression causes some clipping at local extrema in the reconstructed force $\tilde g$, the overall trend of the exact data is captured effectively.
The $L^1-$relative error at terminal time is $\displaystyle \frac{||g(T)-\tg(T)||_1}{||g(T)||_1} = 3.37\times10^{-2}$. The history of state error, $||u(t)-v(t)||_1$, is shown in Figure \ref{fig:linear err}. The result confirms the exponential convergence rate, as expected from convergence analysis.

\subsection{Nonlinear system}
We now examine the performance of the algorithm when applied to the one-dimensional isentropic Euler system for gas dynamics:
\[
\frac{\partial U}{\partial t}+\frac{\partial }{\partial x}F(U) = G,
\]
with
\[
U = 
\left[\begin{array}{c}
\rho\\
\rho u
\end{array}\right], \quad
F(U) = 
\left[\begin{array}{c}
\rho u\\
\rho u^2
\end{array}\right],\quad
G = 
-\frac{\partial}{\partial x}
\left[\begin{array}{c}
0\\
p
\end{array}\right].
\]
Here $\rho$ is the density, $m = \rho u$ is the momentum, and  $p$ is the pressure given by the power law
\[
p = \kappa \rho^\gam, \quad \kappa=1,\hspace{0.2em} \gamma=1.4.
\]
The exact data is generated with the initial condition
\[
\rho_0(x) = 1+0.2\sin(\pi x), \quad u_0(x) = 1,
\]
over the computational domain $[0,4]$. Periodic boundary conditions are applied at both endpoints. Assuming sparse observations are available for the density and the momentum, we expect to recover the true data with the nudged equations,
\[
\left
\{\begin{array}{l}\displaystyle
\frac{\partial V}{\partial t}+\frac{\partial }{\partial x}F(V) = \tG+\mu(I_h(U)-I_h(V)) \,,\\
\\
\displaystyle
\tG = \frac{\partial}{\partial t} I_h(U)+I_h(\frac{\partial}{\partial x}F(\tU)), \quad \tU = I_h(U)+(I-I_h)(V)\,.
\end{array}
\right.
\]
The two components of $V$ are the approximate density $\widehat{\rho}$ and momentum $\widehat{m} = \widehat{\rho u}$. The reconstructed pressure gradient, $\widetilde{\partial_x p}$, is given by the second component of $-\tG$. The nudged system is initialized with 
\[
\widehat{\rho}_0(x) = 1, \quad \widehat{u}_0(x) = 0.5, \quad \tG_0(x) = \left[\begin{array}{c}
0\\
0
\end{array}\right].
\]
The nudging coefficient is set to $\mu = 5$. Numerical discretization is implemented over $N=600$ mesh points.

\begin{figure}[h!]
    \centering
    \begin{subfigure}{0.3\textwidth}
    \includegraphics[scale = 0.22]{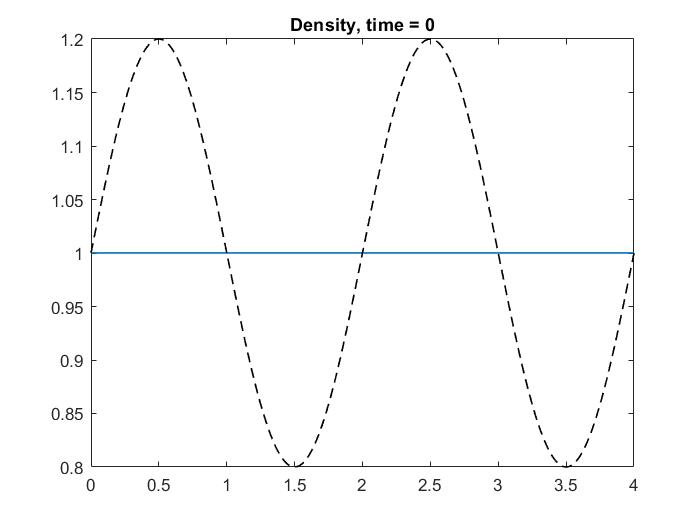}
    \subcaption{$t = 0$}
    \end{subfigure}
    \quad
    \begin{subfigure}{0.3\textwidth}
    \includegraphics[scale = 0.22]{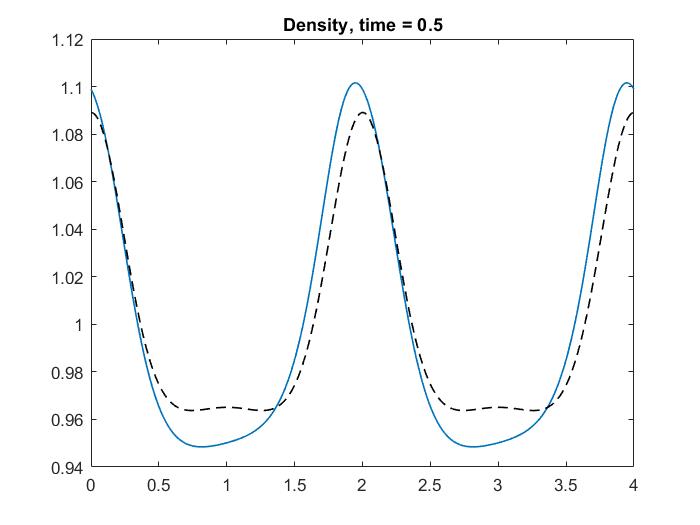}
    \subcaption{$t = 0.5$}
    \end{subfigure}
    \\
    \begin{subfigure}{0.3\textwidth}
    \includegraphics[scale = 0.22]{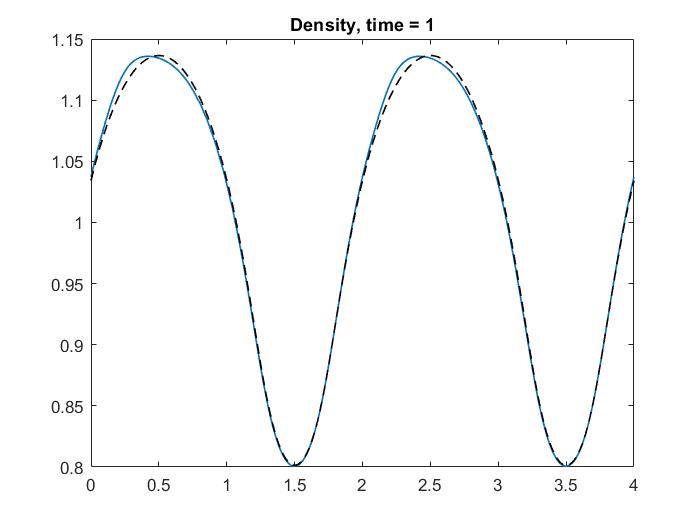}
    \subcaption{$t = 1$}
    \end{subfigure}
    \quad
    \begin{subfigure}{0.3\textwidth}
    \includegraphics[scale = 0.22]{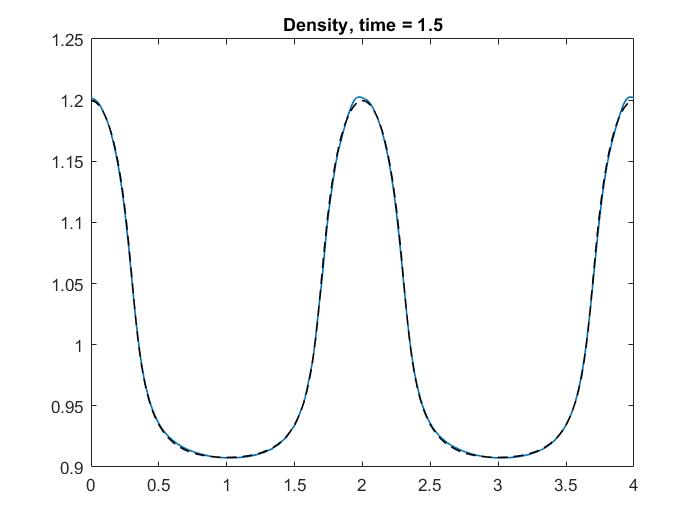}
    \subcaption{$t = 1.5$}
    \end{subfigure}
    \caption{1D isentropic Euler. Density recovery. Nudged solution (blue lines) vs. reference data (black dash lines).}
    \label{fig:isen euler dens}
\end{figure}

\begin{figure}[h!]
    \centering
    \begin{subfigure}{0.3\textwidth}
    \includegraphics[scale = 0.23]{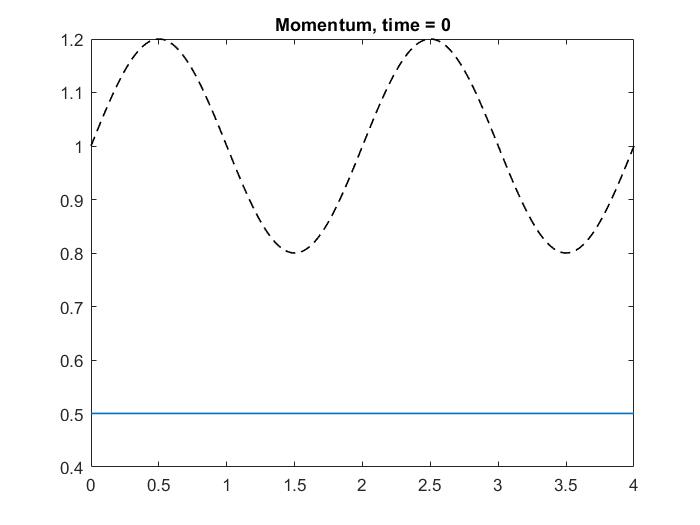}
    \subcaption{$t = 0$}
    \end{subfigure}
    \quad
    \begin{subfigure}{0.3\textwidth}
    \includegraphics[scale = 0.23]{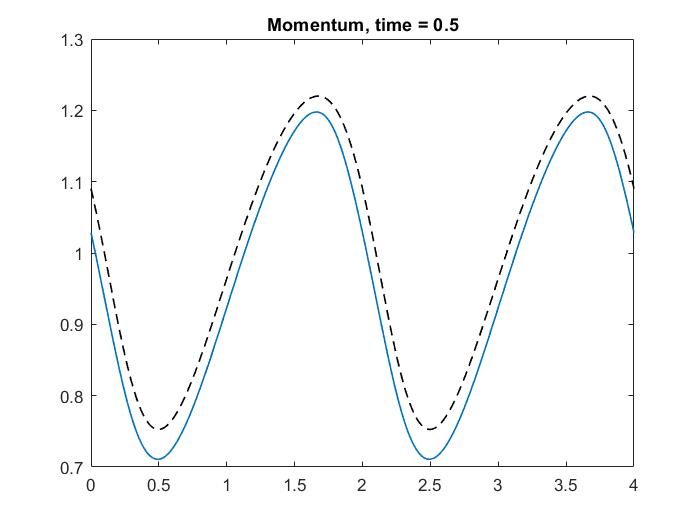}
    \subcaption{$t = 0.5$}
    \end{subfigure}
    \\
    \begin{subfigure}{0.3\textwidth}
    \includegraphics[scale = 0.23]{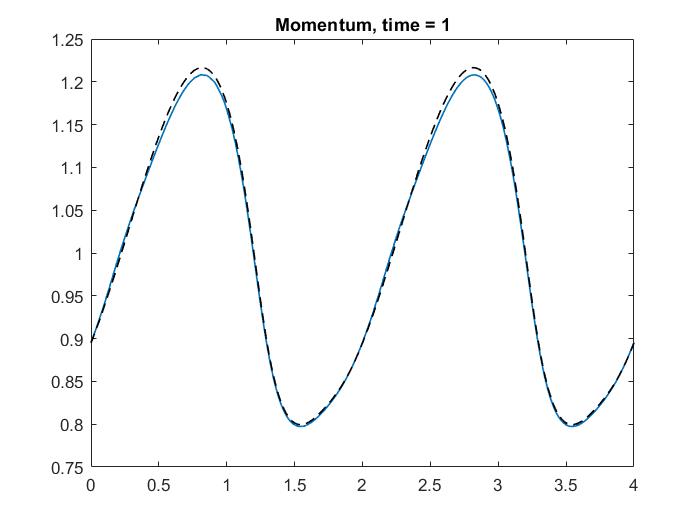}
    \subcaption{$t = 1$}
    \end{subfigure}
    \quad
    \begin{subfigure}{0.3\textwidth}
    \includegraphics[scale = 0.23]{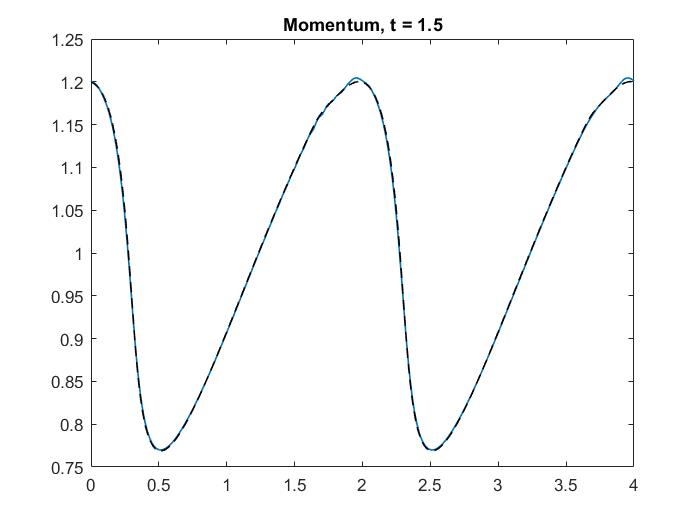}
    \subcaption{$t = 1.5$}
    \end{subfigure}
    \caption{1D isentropic Euler. Momentum recovery. Nudged solution (blue lines) vs. reference data (black dash lines).}
    \label{fig:isen euler mom}
\end{figure}

\begin{figure}[h!]
    \centering
    \begin{subfigure}{0.3\textwidth}
    \includegraphics[scale = 0.23]{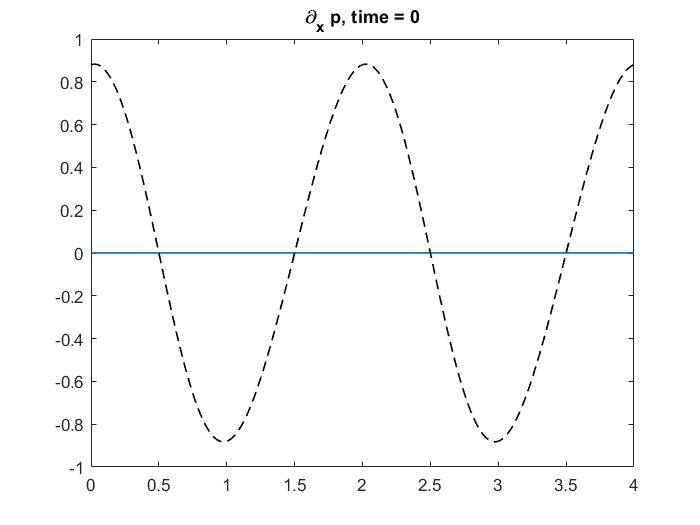}
    \subcaption{$t = 0$}
    \end{subfigure}
    \quad
    \begin{subfigure}{0.3\textwidth}
    \includegraphics[scale = 0.23]{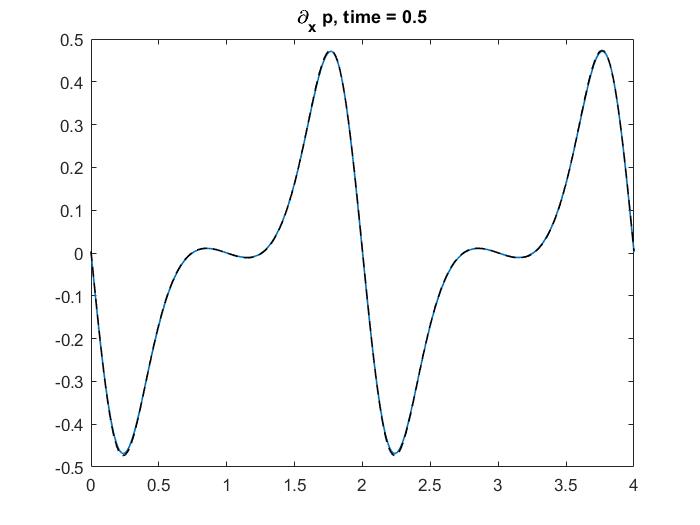}
    \subcaption{$t = 0.5$}
    \end{subfigure}
    \\
    \begin{subfigure}{0.3\textwidth}
    \includegraphics[scale = 0.23]{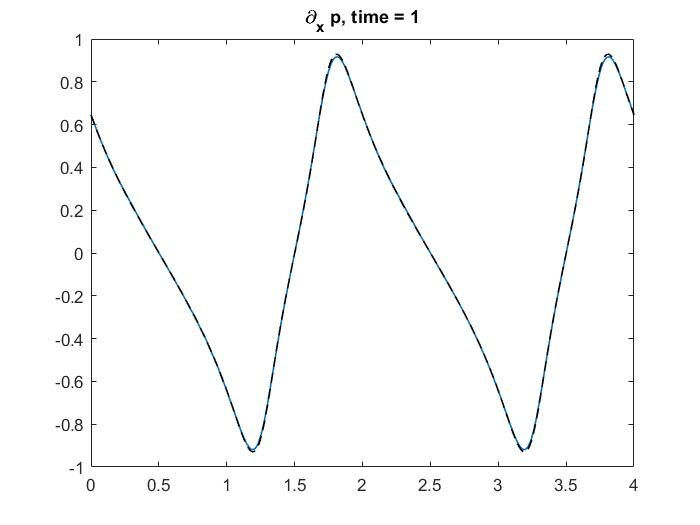}
    \subcaption{$t = 1$}
    \end{subfigure}
    \quad
    \begin{subfigure}{0.3\textwidth}
    \includegraphics[scale = 0.23]{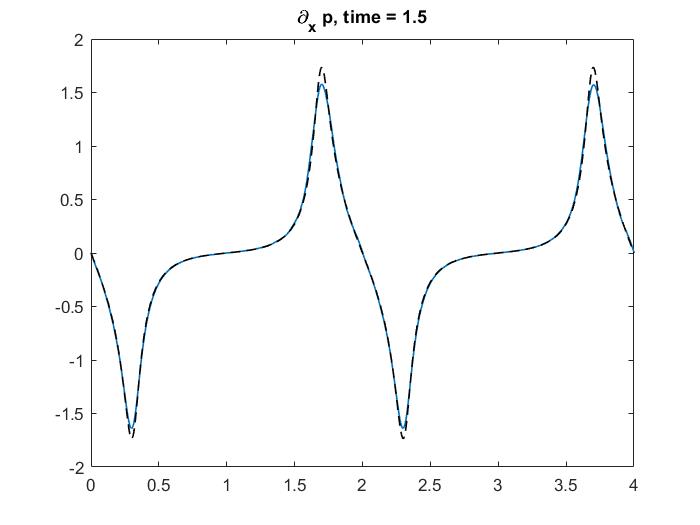}
    \subcaption{$t = 1.5$}
    \end{subfigure}
    \caption{1D isentropic Euler. Pressure gradient recovery. Nudged solution (blue lines) vs. reference data (black dash lines).}
    \label{fig:isen euler pres}
\end{figure}

\begin{figure}[h!]
    \centering
    \begin{subfigure}{0.35\textwidth}
    \includegraphics[scale = 0.24]{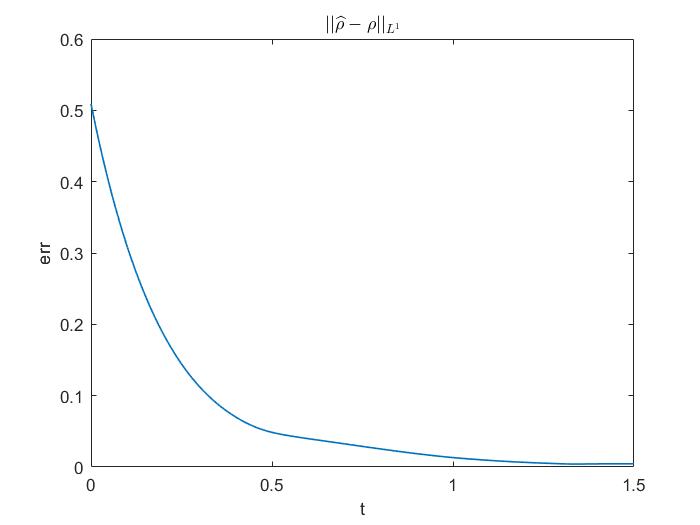}    
    \subcaption{Density error}
    \end{subfigure}
    \quad
    \begin{subfigure}{0.35\textwidth}
    \includegraphics[scale = 0.24]{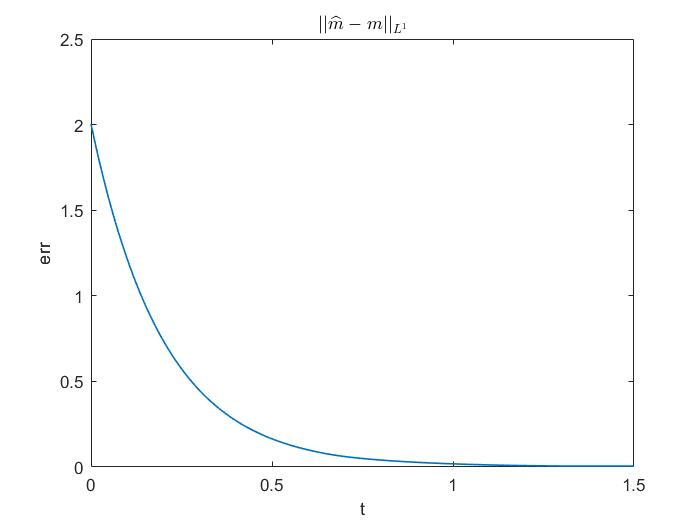}    
    \subcaption{Momentum error}
    \end{subfigure}
    \caption{1D isentropic Euler. History of state errors.}
    \label{fig:isen euler err}
\end{figure}

Figures \ref{fig:isen euler dens} -- \ref{fig:isen euler pres} present the approximations to each state component and the pressure gradient based on $N_{ob}=150$  equidistant observation grids up to $T = 1.5$. The results have shown a stable convergence of computed solutions towards the exact data. The recovered pressure gradient, although slightly clipped at sharp extrema, shows a small terminal relative error of $\displaystyle \frac{||\partial_x p(T)-\widetilde{\partial_x p}(T)||_1}{||\partial_x p(T)||_1} = 3.56\times 10^{-2}$. The state errors of density and momentum are shown in Figure \ref{fig:isen euler err}, which again verify the exponential convergence rates. 

\subsection{Multi-dimensional system}
To validate our algorithm in multi-dimensional problems, we consider the two-dimensional isentropic Euler equations
\[
\begin{split}
\frac{\partial U}{\partial t}+\frac{\partial}{\partial x} F_x(U)+\frac{\partial}{\partial y} F_y(U) = G 
\end{split}
\]
with
\[
U = \left[
\begin{array}{c}
     \rho \\
     \rho u\\
     \rho v
\end{array}\right], \quad
F_x(U) = \left[
\begin{array}{c}
     \rho u\\
     \rho u^2\\
     \rho u v
\end{array}\right], \quad
F_x(U) = \left[
\begin{array}{c}
     \rho v\\
     \rho u v\\
     \rho v^2
\end{array}\right], \quad
G = -\left[
\begin{array}{c}
     0\\
     \partial_x p\\
     \partial_y p
\end{array}\right].
\]
The system is closed by the following constitutive relation:
\[
p = \kappa\rho^\gam, \quad \kappa = 0.5,\hspace{0.2em} \gam = 1.2.
\]
The exact data is generated with the initial conditions
\[
\rho_0(x,y) = 1+0.4\sin(\pi x)\cos(\pi y),\quad u_0(x,y) = 1,\quad v_0(x,y) = 0.5, \quad x,y \in [-1,1]^2.
\]
Periodic boundary conditions are applied at the four boundaries. The nudged solutions are generated with the dynamics
\[
\left\{\begin{array}{l}\displaystyle
\frac{\partial V}{\partial t}+\frac{\partial }{\partial x}F_x(V)+\frac{\partial}{\partial y}F_y(V) = \tG+\mu(I_h(U)-I_h(V))\ ,\\
\\
\displaystyle
\tG = \frac{\partial}{\partial t} I_h(U)+I_h(\frac{\partial}{\partial x}F_x(\tU)+\frac{\partial}{\partial y}F_y(\tU)), \quad \tU = I_h(U)+(I-I_h)(V)\,.
\end{array}\right.
\]
The three components of $V$ are the approximate density $\widehat{\rho}$,  $x-$momentum $\widehat{m}_x = \widehat{\rho u}$, and $y-$momentum $\widehat{m}_y = \widehat{\rho v}$. The reconstructed pressure derivatives, $\widetilde{\partial_x p}$ and $\widetilde{\partial_y p}$, are given by the second and the third components of $-\tG$. We initialize the nudged system with
\[
\widehat{\rho}_0(x,y) = 1, \quad \widehat{u}_0(x,y) = 1,\quad \widehat{v}_0(x,y) = 0.5,\quad \tG = \left[\begin{array}{c}
0\\
0\\
0
\end{array}\right].
\]
The nudging coefficient is set to $\mu = 8$. The kernel regression operator $I_h $ is obtained with bandwidth $\sig = \max\{h^{ob}_x,h^{ob}_y\}$, where $h_x$ and $h_y$ are the mesh sizes of observation grids along $x-$ and $y-$directions. The exact data are observed over the $60\times60$ observation mesh. The nudged equations are discretized with WENO5-SSPRK3 under $180\times180$ computational mesh. 

\begin{figure}[h!]
    \centering
    \begin{subfigure}{0.3\textwidth}
    \includegraphics[scale = 0.6]{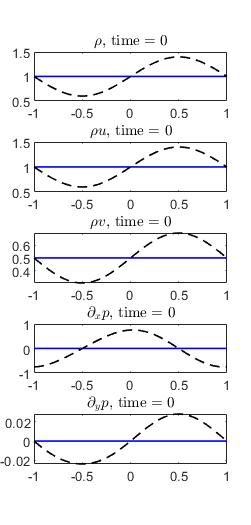}   
     \subcaption{$t = 0$}
    \end{subfigure}
    \quad
    \begin{subfigure}{0.3\textwidth}
    \includegraphics[scale = 0.6]{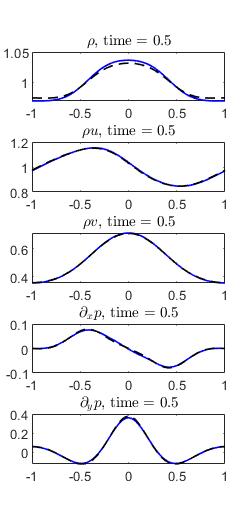}
    \subcaption{$t = 1$}
    \end{subfigure}
    \begin{subfigure}{0.3\textwidth}
    \includegraphics[scale = 0.6]{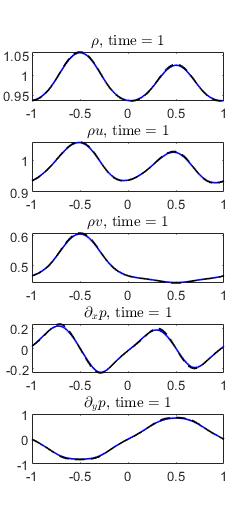}    
    \subcaption{$t = 1.5$}
    \end{subfigure}
    
    \caption{2D Euler. Nudged solutions (blue lines) vs. exact data (black dash lines). }
    \label{fig:2d isen euler}
\end{figure}

\begin{figure}[h!]
    \centering
    \begin{subfigure}{0.3\textwidth}
    \includegraphics[scale = 0.24]{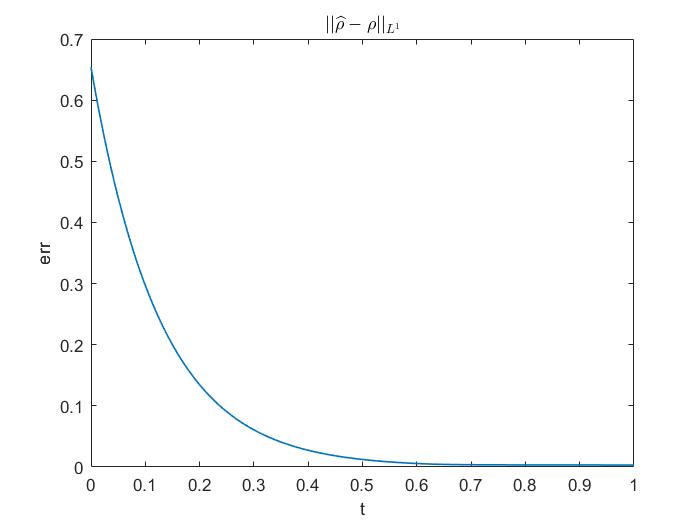}   
    \subcaption{Density error}
    \end{subfigure}
    \quad
    \begin{subfigure}{0.3\textwidth}
    \includegraphics[scale = 0.24]{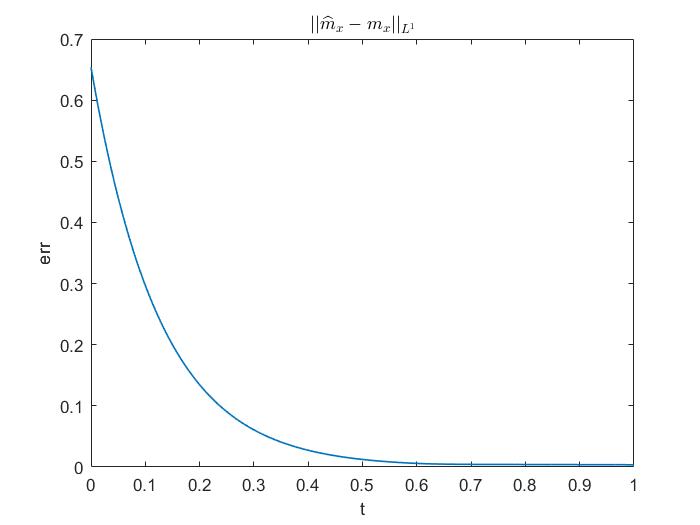}   
    \subcaption{$x-$momentum error}
    \end{subfigure}
    \quad
    \begin{subfigure}{0.3\textwidth}
    \includegraphics[scale = 0.24]{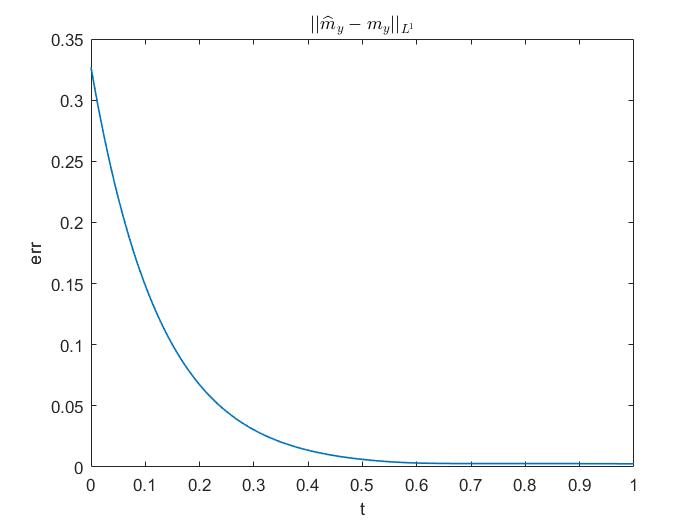}   
    \subcaption{$y-$momentum error}
    \end{subfigure}
    
    \caption{2D Euler. History of state errors.  }
    \label{fig:2d isen euler err}
\end{figure}

Figure \ref{fig:2d isen euler} presents the computed solutions along $y = 0$ up to $T = 1$. The results verify the quick convergence of the approximate density and momentum towards the reference solutions. The computed pressure gradients closely characterize the variation in the exact data. The terminal relative errors are $\displaystyle\frac{||\widetilde{\partial_x p}(T)-\partial_x p(T)||_1}{||\partial_x p(T)||_1} = 1.80\times10^{-2}$ and $\displaystyle\frac{||\widetilde{\partial_y p}(T)-\partial_y p(T)||_1}{||\partial_y p(T)||_1} = 1.79\times10^{-2}$. The history of state errors, as shown in Figure \ref{fig:2d isen euler err}, verifies the exponential convergence in multi-dimensions.

\vspace{0.3in}

\subsection{Noisy measurement and sparse observation}
\label{subsec:noisy}

In practical implementations of data assimilation algorithms, challenges often arise due to noisy and/or sparse observations, which can degrade the accuracy of data recovery. In this section, we primarily investigate the performance of our algorithm under noisy observational data and varying observation density.

\begin{figure}[h!]
    \centering
    \begin{subfigure}{0.31\textwidth}
    \includegraphics[scale = 0.18]{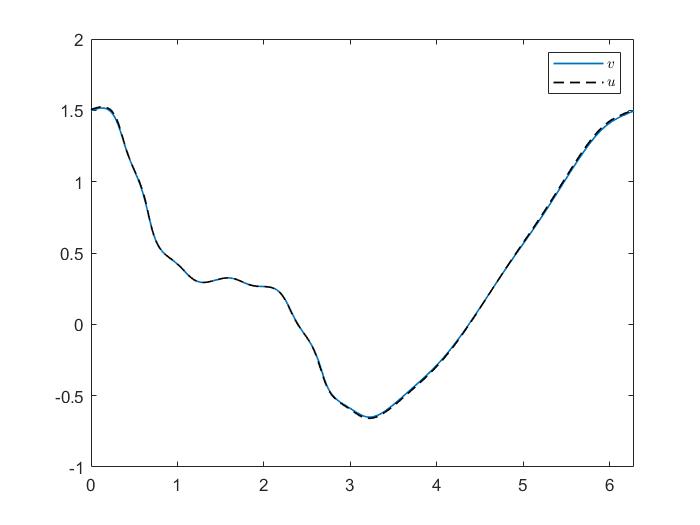}
    \subcaption{$\sig = h$}
    \end{subfigure}
    \quad 
    \begin{subfigure}{0.31\textwidth}
    \includegraphics[scale = 0.18]{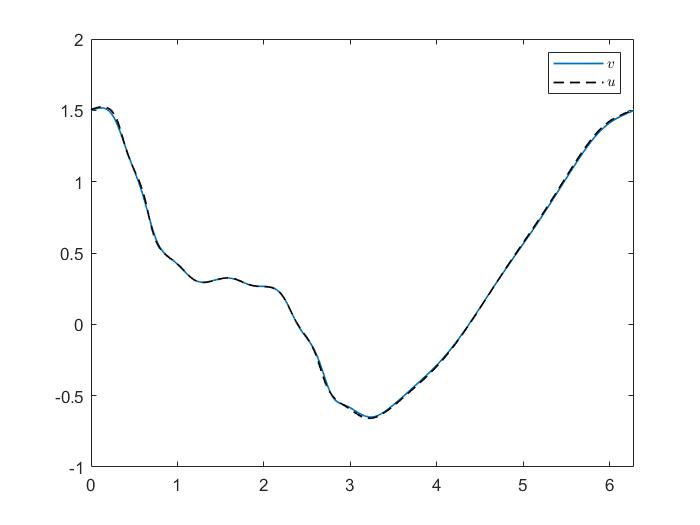}
    \subcaption{$\sig = 2h$}
    \end{subfigure}
    \begin{subfigure}{0.31\textwidth}
    \includegraphics[scale = 0.18]{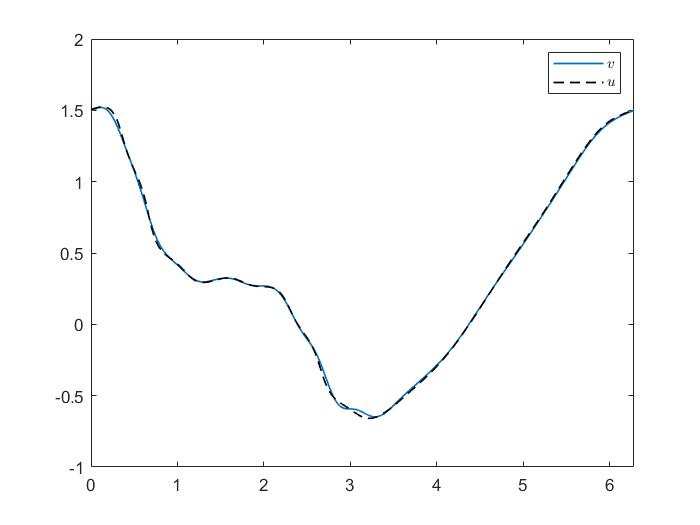}
    \subcaption{$\sig = 3h$}
    \end{subfigure}
    \\
    \begin{subfigure}{0.31\textwidth}
    \includegraphics[scale = 0.18]{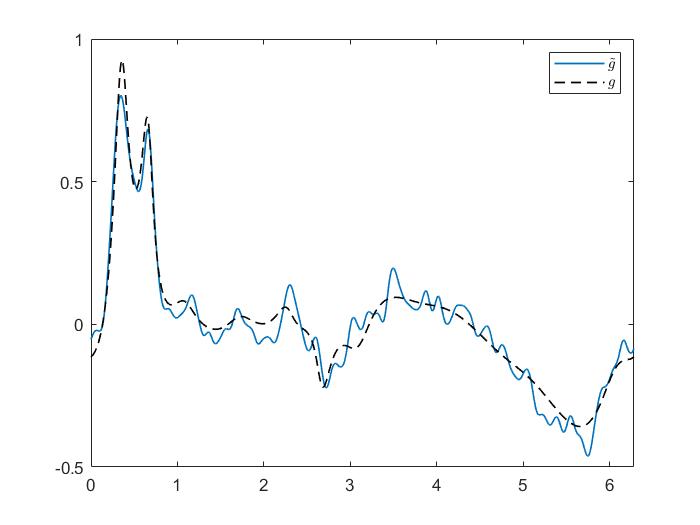}
    \subcaption{$\sig = h$}
    \end{subfigure}
    \quad 
    \begin{subfigure}{0.31\textwidth}
    \includegraphics[scale = 0.18]{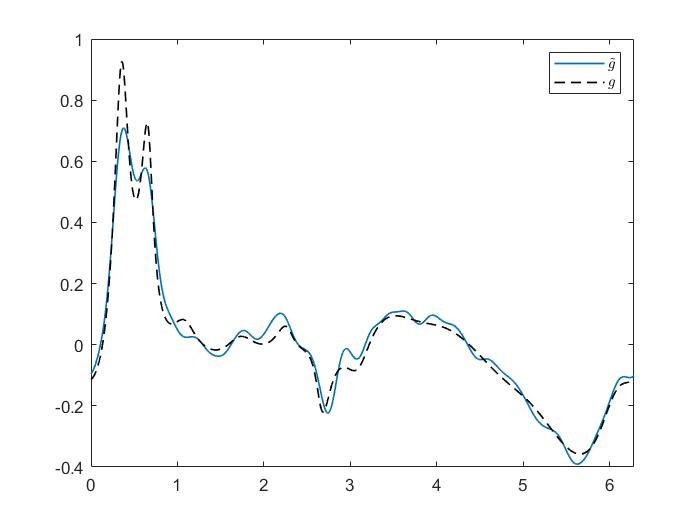}
    \subcaption{$\sig = 2h$}
    \end{subfigure}
    \begin{subfigure}{0.31\textwidth}
    \includegraphics[scale = 0.18]{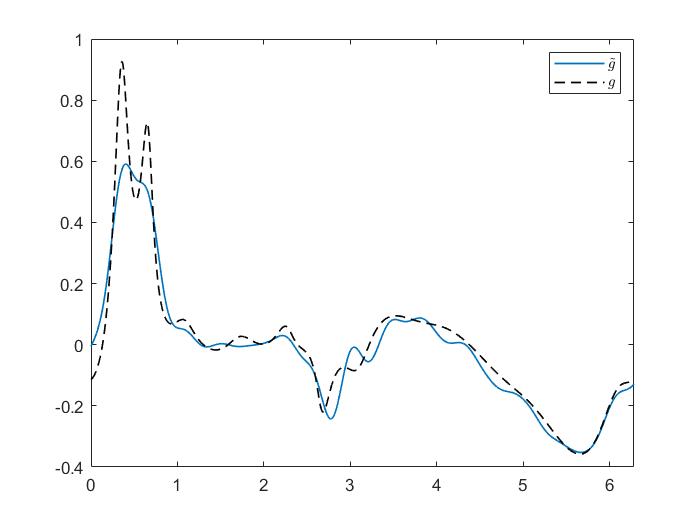}
    \subcaption{$\sig = 3h$}
    \end{subfigure}
    \caption{Scalar dynamics with noisy observation. Nudged solutions (blue lines) vs. reference solutions (black dash lines) at $t = 1.5$ under varying regression bandwidth. The first row presents the results for state recovery, and the second row presents the results for force recovery.}
    \label{fig:linear with noise}
\end{figure}

For convenience of demonstration, we consider the scalar dynamics \eqref{eq:scalar dynamics} and the corresponding nudging system \eqref{eq:scalar nudging}. At first, we examine the performance of our data assimilation algorithm under noisy observations, namely, the feedback control term is constructed with $I_h(u^\eps(t))$ where
\[
u^\eps(x_j^{ob},t) = u(x^{ob}_j,t)+\eps \mathcal{N}_{j,t},
\]
with $\{\mathcal{N}_{j,t}\}$ independently and identically distributed with respect to the standard normal distribution. 
Following the computational settings in Section \ref{sec:scalar test}, we still take the nudging coefficient $\mu = 3$ and assume that the measurement of $u$ is conducted over $N_{ob} = 150$ grids, i.e., the mesh size of observation grids is $h = 2\pi/150$. We take the magnitude of noise to $\eps = 10^{-4}$. 
As demonstrated in Section \ref{sec:kernel regression}, a key advantage of using kernel regression as an observable interpolant lies in its convenience and robustness in processing noisy data. To illustrate this, we display the computational results in Figure \ref{fig:linear with noise}. It is observed that,  while the recovery of the solution state $u$ remains relatively unaffected by small levels of noise, the recovery of the unknown forcing term $g = -\partial_x p$ is significantly more sensitive to inaccuracies in the observations. By increasing the bandwidth of the regression kernel, the recovered data can be effectively smoothed, although at the cost of more dissipative local extrema. In practical applications, the value of $\sigma$ can be tuned to strike a balance between preserving sharp features and achieving desired smoothness in the reconstructed solutions.

\begin{figure}[h!]
    \centering
    \begin{subfigure}{0.3\textwidth}
    \includegraphics[scale = 0.18]{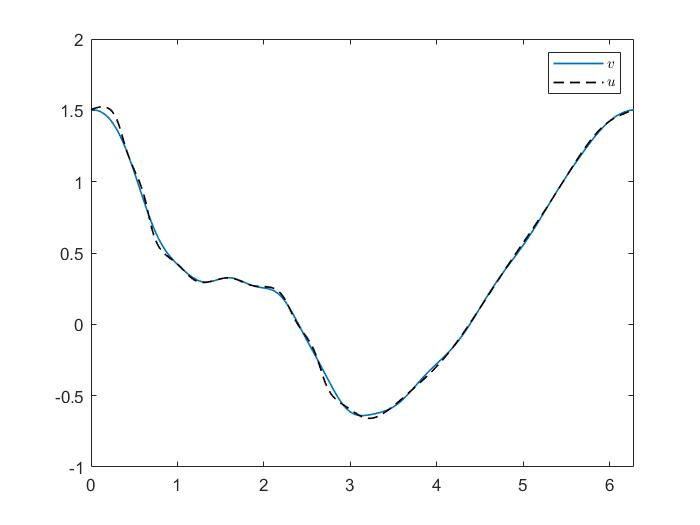}
    \subcaption{$N_{ob}=30$}
    \end{subfigure}
    \quad 
    \begin{subfigure}{0.3\textwidth}
    \includegraphics[scale = 0.18]{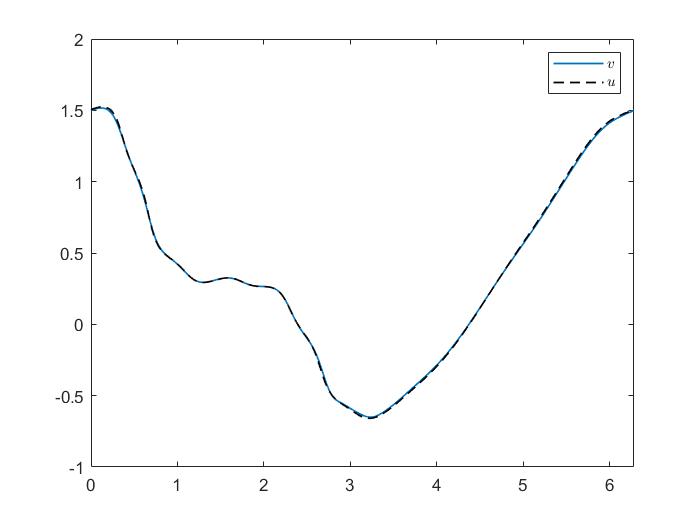}
    \subcaption{$N_{ob}=90$}
    \end{subfigure}
    \quad 
    \begin{subfigure}{0.3\textwidth}
    \includegraphics[scale = 0.18]{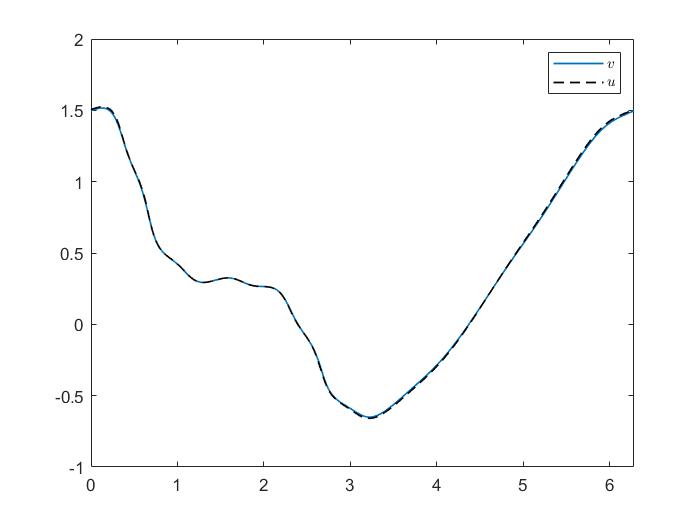}
    \subcaption{$N_{ob}=150$}
    \end{subfigure}
    \\
    \begin{subfigure}{0.3\textwidth}
    \includegraphics[scale = 0.18]{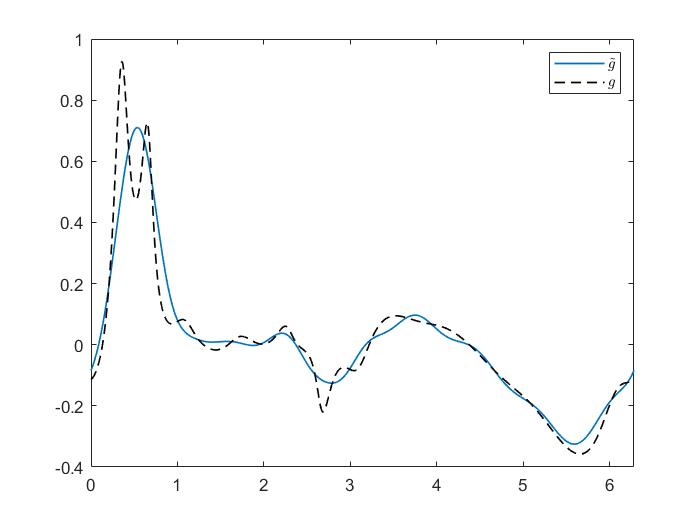}
    \subcaption{$N_{ob}=30$}
    \end{subfigure}
    \quad 
    \begin{subfigure}{0.3\textwidth}
    \includegraphics[scale = 0.18]{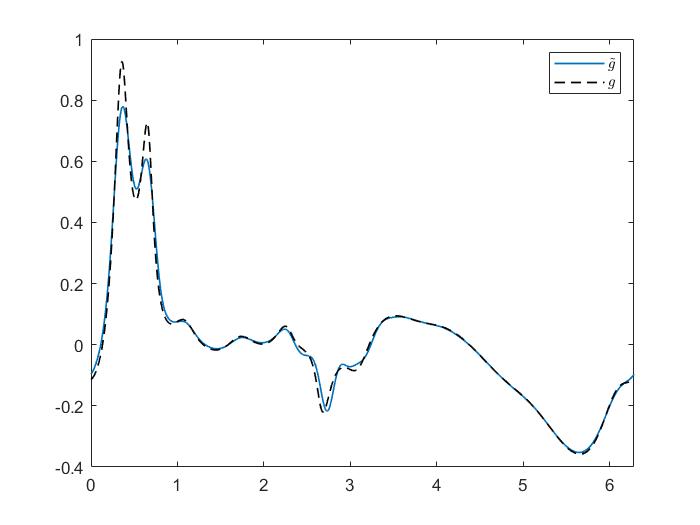}
    \subcaption{$N_{ob}=90$}
    \end{subfigure}
    \quad 
    \begin{subfigure}{0.3\textwidth}
    \includegraphics[scale = 0.18]{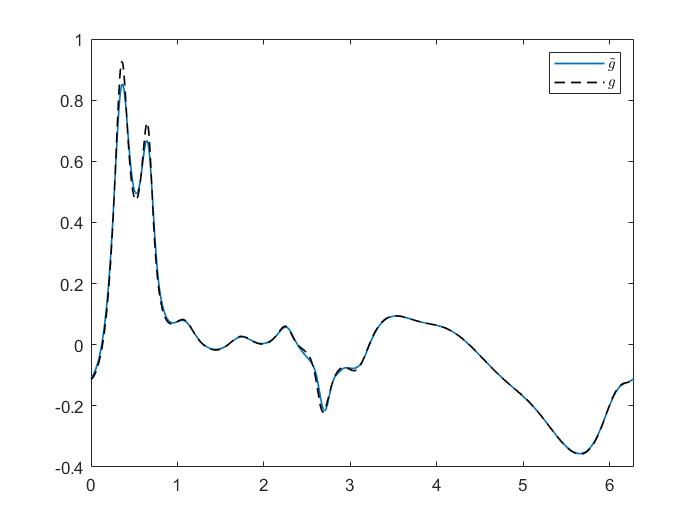}
    \subcaption{$N_{ob}=150$}
    \end{subfigure}
    \caption{Scalar dynamics with varying observation density. Nudged solutions (blue lines) vs. reference solutions (black dash lines) at $t = 1.5$ under varying observation density. The first row presents the results for state recovery, and the second row presents the results for force recovery.}
    \label{fig:linear sparse observ}
\end{figure}

We also evaluate the performance of our algorithm under sparse observational data. In each test case, the bandwidth of the regression kernel is set to $\sig = h = 2\pi/N_{ob}$, with $N_{ob}$ being the number of observation grids. As shown in Figure \ref{fig:linear sparse observ}, it can be found that reducing the number of observation grids leads to less clear details in the recovered solution. Nevertheless, the overall variation of the solutions can still be reasonably captured, even when the observations are very sparse. In practical scenarios where the number of observation points (sensors) is limited, local refinement of the observation grids in the regions with complex interactions may further improve the reconstruction accuracy. A detailed investigation of such adaptive strategies is left for our future work.

\section{Moment recovery with incomplete observations}
\label{sec:moment-recovery}

The previous discussion has established the validity of algorithm \eqref{eq:nudge} when observational data for all state components are available, a situation we referred to as \emph{complete} observation. However, in some applications, only some of the state components may be observable. In such cases, where data for some state components are missing, we refer to the observation as \emph{incomplete}. This scenario is particularly relevant in kinetic theory, where one aims to derive a moment system as a reduced-order model. In this section, we demonstrate that the method proposed in Section 2 can be extended to handle incomplete observations, providing an efficient mechanism for moment recovery.

\subsection{Problem set up }
As a concrete example, we consider the one-dimensional radiative transfer equation (RTE):
\begin{equation}\label{eq:rte}
    \partial_t f + v\partial_x f = \sig_s(\frac{1}{2}\int^1_{-1} f \rd v-f)-\sig_a f, \quad x\in [0,1], \hspace{0.5em} v\in[-1,1].
\end{equation}
Here $f(x,v,t)$ is the specific intensity of radiation, $v\in[-1,1]$ is the cosine of the angle between the photon velocity and the $x$-axis. $\sig_s(x)\geq 0$ and $\sig_a(x)\geq 0$ are the scattering and absorption coefficients, respectively. 

Let $P_k(v)$ be the Legendre polynomial of $k$ degrees. Integrating the equation \eqref{eq:rte} against $P_k(v)$ over $v\in[-1,1]$ and using Bonnet's recursion formula, one obtains the moment equations up to the $N-$th order:

\begin{equation}\label{eq:rte moment}
\begin{split}
\partial_t m_0 + \partial_x m_1 &= -\sig_a m_0, \\
\partial_t m_1 +\frac{1}{3}\partial_x m_0+\frac{2}{3}\partial_x m_2 &= -(\sig_s+\sig_a) m_1,\\
&\cdots  \\
\partial_t m_N +\frac{N}{2N+1}\partial_x m_{N-1}+\frac{N+1}{2N+1}\partial_x m_{N+1} &= -(\sig_a+\sig_s)m_N,
\end{split}
\end{equation}
where the $k-$th moment is defined as
\[
m_k(x,t) = \frac{1}{2}\int^1_{-1}f(x,v,t)P_k(v)dv.
\]
The system \eqref{eq:rte moment} is deemed a model reduction for \eqref{eq:rte}, but it is not closed as the evolution of $m_N$ depends on the unknown moment $m_{N+1}$. One approach to closure involves using a machine learning-based moment closure \cite{huang2022machineI,huang2023machineII,huang2023machineIII}. Given sufficient data of moments $m_1$ to $m_N$, which can be obtained by solving \eqref{eq:rte}, one attempts to learn a relationship such as:
\begin{equation}\label{eq:ml closure}
\partial_x m_{N+1} = \sum^{N}_{k=0}\network_k(m_0,m_1,\cdots,m_N)\partial_x m_k.
\end{equation}

However, direct simulation of the kinetic equation \eqref{eq:rte} can be expensive, and instead we aim to reconstruct $\partial_x m_{N+1}$ from the observed data of lower-order moments. To this end, we write the moment equations in the condensed form:
\begin{equation}\label{eq:rte moment system}
\begin{split}
    &\frac{\partial \bm_{N}}{\partial t}+A_N\frac{\partial\bm_{N}}{\partial x} = \bg_{N+1}+S_N\bm_{N}, \\
    &\bg_{N+1} = -\frac{N+1}{2N+1}\partial_x m_{N+1} \be_{N+1},
\end{split}
\end{equation}
with $\bm_{N} = (m_0,m_1,\cdots,m_N)^\top$, $\be_{N+1} = (0,\cdots,0,1)^\top$, and the coefficient matrices:
\begin{align} \label{eqn:A}
A_N = \left[\begin{array}{cccccc}\
    0 & 1 & 0 & 0 & \cdots & 0 \\ 
    \frac{1}{3} & 0 & \frac{2}{3} & 0 &\cdots & 0 \\
    0 & \frac{2}{5} & 0 & \frac{3}{5} & \cdots & 0\\
    \vdots & \vdots & \vdots & \ddots & \vdots & \vdots \\
    0 & 0 & \cdots & \frac{N-1}{2N-1} & 0 & \frac{N}{2N-1}\\ 
    0 & 0 & \cdots & 0 & \frac{N}{2N+1} & 0
\end{array}\right]\in \R^{(N+1)\times(N+1)},
\end{align}

\begin{align} \label{eqn:S}
S_N = \text{diag}(-\sig_a,-(\sig_s+\sig_a),\cdots,-(\sig_s+\sig_a))\in \R^{(N+1)\times(N+1)}.
\end{align}

Ideally, if all components of $\bm_{N}$ are measurable at sparse observation points, one could asymptotically reconstruct the continuous solution fields $\bm_{N}(x,t)$ and $\partial_x m_{N+1}(x,t)$ over time, following the methodology described in \eqref{eq:nudge}. These reconstructed data could then be utilized to train neural networks to model the relationship in \eqref{eq:ml closure}. However, obtaining data for all moments is impractical in real-world applications. Indeed, higher-order moments beyond the first two—angle-averaged intensity $m_0$ and Eddington flux $m_1$—often lack direct physical interpretation and are not readily available for measurement. Nonetheless, \cite{huang2022machineI} highlights the necessity of using a sufficient number of moments $(N \geq 5)$  to ensure accurate results in both optically-thin and optically-thick regimes. Therefore, it is crucial to obtain data for $m_2$, $m_3$, $\ldots$, $m_{N+1}$. 
To avoid the challenges associated with microscopic or full measurements of the particle distribution $f$ in complex applications, we propose using the continuous data assimilation mechanism described earlier.

\subsection{Moment recovery via continuous data assimilation}
Back to the moment system \eqref{eq:rte moment} (or \eqref{eq:rte moment system}), we distinguish the notations of the ground truth data, $m_k(x,t)$, and the numerical solutions, $\hm_k(x,t)$. Our purpose is to recover the continuous data fields of all moments, $\bm_{N+1} = (m_0,m_1,\cdots,m_{N+1})^\top$, while using the expensive observations of higher-order moments to a minimal extent. To this end, we assume that the full measurements of $f$, and thereby complete observations of all moments of interest, are {\it only conducted at the boundaries}. At the interior of the solution domain, we employ incomplete observations with the macroscopically observable moments up to the $n-$th order, $\bm_n = (m_0,m_1,\cdots,m_n)^\top$, $n<N$. Moreover, the gradients $\partial_x m_{n+1}, \partial_x m_{n+2}, \cdots, \partial_x m_{N+1}$ are assumed to be completely observable with respect to the interpolant operator $I_h$, i.e., 
\[
I_h(\partial_x m_k) = \partial_x m_k, \quad k = n+1, n+2, \cdots N.
\] \
We propose the following moment recovery algorithm for the one-dimensional RTE.

\medskip
\medskip
\paragraph{Step 1.} Reconstruct the moments $\widehat{\bm}_n(x,t) = (\hm_0(x,t), \hm_1(x,t),\cdots,\hm_n(x,t))^\top$ by solving the nudged system with observations $\bm_n$:
\[
\left\{\begin{array}{l}
\displaystyle
\frac{\partial \widehat{\bm}_n}{\partial t}+A_n\frac{\partial \widehat{\bm}_n}{\partial x} = \widetilde{\bg}_{n+1}+S_n\widehat{\bm}_n+\mu(I_h(\bm_n)-I_h(\widehat{\bm}_n)),  \\
\\
\displaystyle
\widetilde{\bg}_{n+1} = \frac{\partial}{\partial t}I_h(\bm_n)+I_h(A_n\frac{\partial \widetilde{\bm}_n}{\partial x}-S_n\widetilde{\bm}_n), \quad \widetilde{\bm}_n = I_h(\bm_n)+(I-I_h)(\widehat{\bm}_n)\,.
\end{array}\right.
\]

\paragraph{Step 2.} Obtain the reconstructed gradient, $\widetilde{\partial_x m}_{n+1}$, from the last component of $\displaystyle -\frac{2n+1}{n+1}\widetilde{\bg}_{n+1}$. Then given the boundary measurement $m_{n+1}(0,t)$ and $m_{n+1}(1,t)$, the approximation to $m_{n+1}$ is then obtained with integration:
\[
\hm_{n+1}(x,t) = m_{n+1}(0,t)+\int^x_0 (\widetilde{\partial_x m}_{n+1}+c_{n+1})dx\,.
\]
Here $c_{n+1}$ is the gradient correction 
\[
c_{n+1} = -\int^1_0 \widetilde{\partial_x m}_{n+1} dx+ m_{n+1}(1,t)- m_{n+1}(0,t),
\]
such that the boundary conditions, $\hm_{n+1}(0,t) = m_{n+1}(0,t)$ and $\hm_{n+1}(1,t) = m_{n+1}(1,t)$, are ensured.  

\paragraph{Step 3.} For $k = n+1, n+2,\cdots,N$, suppose the nudged solutions $\hm_0(x,t), \cdots, \hm_k(x,t)$ have been obtained, to get $\hm_{k+1}(x,t)$, we  first compute its approximate gradient with
\[
\widetilde{\partial_x m}_{k+1} = -\frac{2k+1}{k+1}I_h(\partial_t \hm_{k}+\frac{k}{2k+1}\partial_x \hm_{k-1}+(\sig_a+\sig_s)\hm_k).
\]
Then we recover $\hm_{k+1}$ through integration
\[
\hm_{k+1}(x,t) = m_{k+1}(0,t)+\int^x_0 (\widetilde{\partial_x m}_{k+1}+c_{k+1}) dx
\]
with 
\[
c_{k+1} = -\int^1_0 \widetilde{\partial_x m}_{k+1} dx+ m_{k+1}(1,t)- m_{k+1}(0,t).
\]
The reconstructed moment satisfies the given boundary measurements: $\hm_{k+1}(0,t) = m_{k+1}(0,t)$ and $\hm_{k+1}(1,t) = m_{k+1}(1,t)$.

We now comment on the availability of the boundary data $m_k(0,t)$ and $m_k(1,t)$. While obtaining measurements of $f$
throughout the entire interior domain can be impractical, measurements at the boundaries are more feasible. Specifically, one can impose a user-specified incoming boundary condition, i.e., $f(0, v>0)$ and $f(1,v<0)$, and then perform velocity-resolved boundary measurements, i.e., $f(0, v<0)$ and $f(1,v>0)$. This provides complete information about $f$ at the boundaries, allowing for the calculation of moments: 
\[
m_k(0) = \frac{1}{2}\int^1_{-1} f(x=0,v) P_k(v)dv, \quad m_k(1) = \frac{1}{2}\int^1_{-1}f(x=1,v) P_k(v)dv.
\]
This input-output data pair is common in inverse problem setups, and their relationship is often referred to as the Albedo operator. For more details, see \cite{Bal_2009}.


It is worth mentioning that the above moment recovery procedure can be naturally extended to general 1-D moment systems that admit the structure:
\[
\partial_t m_k +\sum^{k+1}_{i=0}a_{ki}(m_0,m_1,\cdots,m_k)\partial_x m_{i} = S_k(m_0,m_1,\cdots,m_k), \quad a_{k,k+1} \neq 0, \hspace{0.5em} k = 0,1,2,\cdots.
\]

For instance, we briefly mention the one-dimensional Boltzmann equation with BGK collision \cite{bhatnagar1954model}: 
\[
\frac{\partial f}{\partial t}+v\frac{\partial f}{\partial x} = \frac{1}{\tau}(f_{\mathcal{M}}-f), \quad x\in [0,1], \quad v\in\R\,.
\]
Here $\tau>0$ is the relaxation time, and  $\displaystyle f_{\mathcal{M}} = \frac{\rho}{\sqrt{2\pi \th}}\exp(-\frac{|v-u|^2}{2\th})$ is the Maxwellian.
We define the moments as $\mathbf{\mathfrak{m}}_N = (\rho,u,\th,f_3,f_4,\cdots,f_{N})^\top\in\R^{N+1}$, 
where 
\[
\rho(x,t) = \int_{\R} f(x,v,t) dv, \quad \rho u(x,t) = \int_{\R} vf(x,v,t) dv,\quad \frac{1}{2}\rho \th(x,t) = \int_\R\frac{(v-u)^2}{2}f(x,v,t)dv\,,
\]
and $f_k$ is the coefficient in the Hermite expansion \cite{cai2010numerical,cai2012numerical}:
\[
f(x,v,t) = \sum^{\infty}_{k = 0} f_k(x,t)\mathcal{H}_{\th(x,t),k}(\frac{v-u(x,t)}{\sqrt{\th(x,t)}}), \quad \mathcal{H}_{\th,k}(v) = \frac{1}{\sqrt{2\pi}}\th^{-\frac{k+1}{2}} He_k(v)\exp(-\frac{v^2}{2}),
\]
where $He_k$ being the $k-$th Hermite polynomial. In particular, it can be shown that $f_0 = \rho$ and $f_1 = f_2 = 0$. Inspired by \cite{cai2013globally}, the corresponding moments system are given by 
\[
\frac{\partial \mathbf{\mathfrak{m}}_{N}}{\partial t} +A_N\frac{\partial \mathbf{\mathfrak{m}}_{N}}{\partial x}  = -(N+1)\frac{\partial f_{N+1}}{\partial x} \be_{N+1}-\frac{1}{\tau}P_N \mathbf{\mathfrak{m}}_{N},
\]
where
\begin{small}
\[
A_N = \left[\begin{array}{cccccccccc}
    u & \rho & 0 & 0 & 0& \cdots & \cdots &\cdots& \cdots & 0 \\ 
    \frac{\th}{\rho} & u & 1 & 0 & 0& \cdots & \cdots& \cdots& \cdots & 0\\[4pt]
    0 & 2\th & u & \frac{6}{\rho} & 0& \cdots & \cdots & \cdots&\cdots & 0\\[4pt]
    0& 4f_3 & \frac{\rho\th}{2}& u & 4 & 0 & \cdots & \cdots & \cdots & 0\\[4pt]
    -\frac{\th f_3}{\rho}& 5f_4 &\frac{3f_3}{2}& \th & u& 5 &0  & \cdots &\cdots & 0\\[4pt]
    \vdots & \vdots & \vdots & \vdots &\cdots & \cdots&\cdots &\cdots&\cdots &\cdots\\[4pt]
    -\frac{\th f_{N-2}}{\rho} & Nf_{N-1} & \frac{1}{2}[(N-2)f_{N-2}+\th f_{N-4}] & -\frac{3f_{N-3}}{\rho} &0&  \cdots& 0 &\th&u&N\\[4pt]
    -\frac{\th f_{N-1}}{\rho} & (N+1)f_{N} & \frac{1}{2}[(N-1)f_{N-1}+\th f_{N-3}] & -\frac{3f_{N-2}}{\rho} &0& \cdots&\cdots& 0 &\th& u
\end{array}\right],
\]
\end{small}
and
\[
P_N = \text{diag}(0,0,0,1,\cdots,1).
\]
Despite the non-constant coefficients, $A_N$ remains an unreduced lower Hessenberg matrix for a physical density, $\rho>0$, thereby allowing the recovery of higher moments through integration.


\subsection{Numerical experiments II -- moment recovery}

\subsubsection{One-dimensional RTE}
We first apply our moment recovery algorithm to the moment system \eqref{eq:rte moment} derived from the one-dimensional RTE. The reference solution is generated by solving the equation \eqref{eq:rte} using the WENO5-SSPRK3 with  initial data 
\[
f_0(x,v) = 0.5+\sum^5_{k=1}\frac{1}{k^2}\sin(2k\pi x), \quad x\in[0,1].
\]
Periodic boundary conditions are applied at both endpoints. The number of grid points in space is set to $N = 600$. The velocity space is discretized with 15 Gauss-Legendre quadrature points. The time step is taken to $\Dt = 0.5\Dx$ and the terminal time is $t = 1$.

We consider the data under different absorption coefficients $\sig_a$ and $\sig_s$. In all test cases, the nudged equations are discretized with WENO5-SSPRK3 under $N = 300$ mesh points. The nudging coefficient is set to $\mu = 6$. Complete observations are assumed to be available at both boundaries. We intend to recover the higher-order moments based on the observed data of $m_0$ and $m_1$  over $N_{ob} = 60$ uniform observation grids at the interior. The approximate moments are initialized with $\hm_0 = 0.5$, $\hm_2=0$, $\hm_3=0$, $\hm_4=0$, $\hm_5=0$.

\begin{figure}[h!]
    \centering
    \includegraphics[scale = 0.53]{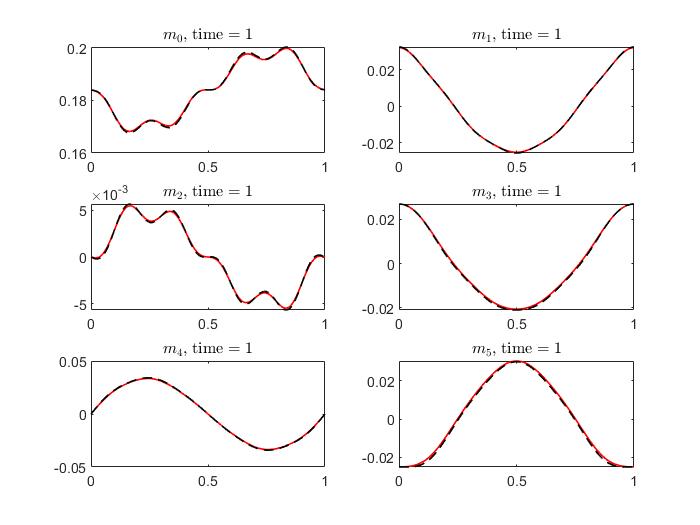}
    \caption{Moment recovery for 1D RTE. Nudged solutions (red lines) vs. reference data (black dash lines). $\sig_a = \sig_s = 1$. $t = 1$.}
    \label{fig:rte1d opt thin}
\end{figure}

\begin{figure}[h!]
    \centering
    \includegraphics[scale = 0.53]{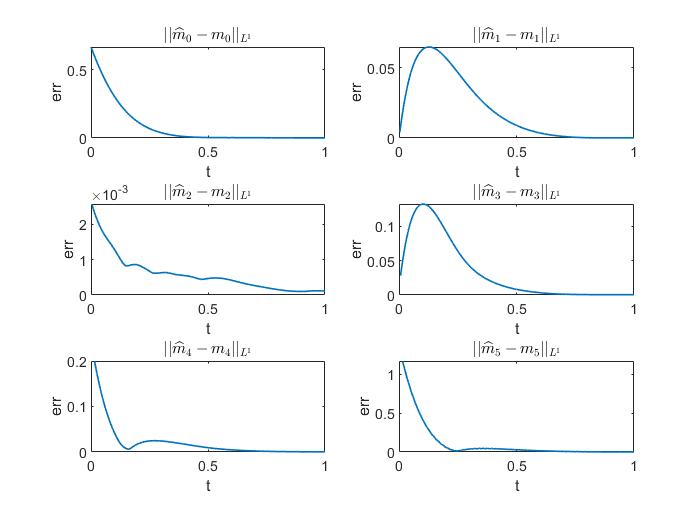}
    \caption{1D RTE. History of moment errors. $\sig_a = \sig_s = 1$.}
    \label{fig:rte1d opt thin err}
\end{figure}

\paragraph{\textit{Test 1 ($\sig_a = \sig_s = 1$).}} We recover the moments up to the fifth order in the optically-thin regime with small absorption coefficients $\sig_a = \sig_s = 1$. The nudged solutions $\hm_k$ at terminal time are displayed in Figure \ref{fig:rte1d opt thin}. It is observed that all the recovered moments match the exact data accurately. 
The history of moment errors, as shown in Figure \ref{fig:rte1d opt thin err}, verifies the convergence of the algorithm, where the relatively large error at early times is primarily due to the intentionally inaccurate initialization used in our test cases. This setup was chosen to demonstrate the robustness of the nudging-based feedback control, particularly its ability to recover accurate solutions even from poor initial guesses.

\begin{figure}[h!]
    \centering
    \includegraphics[scale = 0.53]{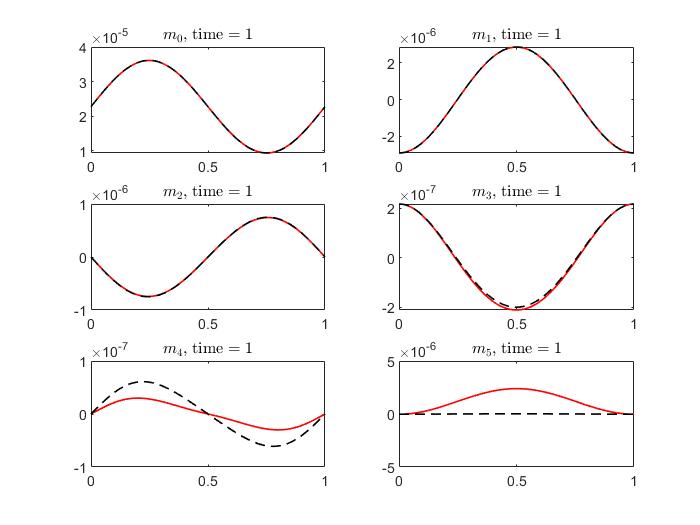}
    \caption{Moment recovery for 1D RTE. Nudged solutions (red lines) vs. reference data (black dash lines). $\sig_a = \sig_s = 10$. $t = 1$.}
    \label{fig:rte1d interm}
\end{figure}

\begin{figure}[h!]
    \centering
    \includegraphics[scale = 0.53]{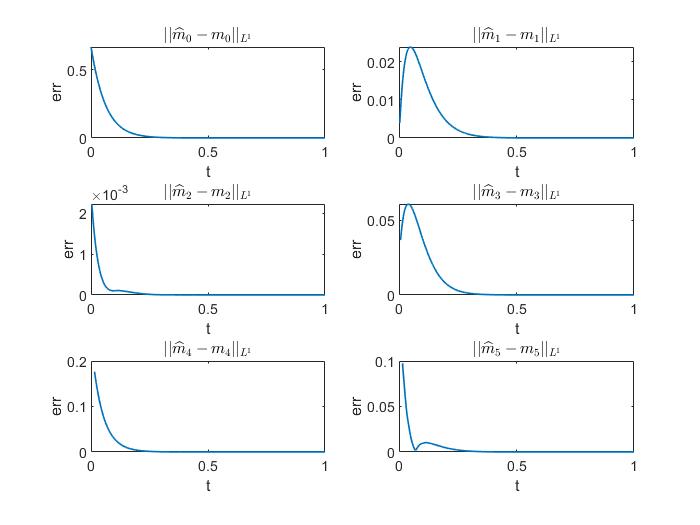}
    \caption{1D RTE. History of moment errors. $\sig_a = \sig_s = 10$.}
    \label{fig:rte1d interm err}
\end{figure}

\paragraph{\textit{Test 2 ($\sig_a = \sig_s = 10$).}} Now we consider the data under a larger optical depth. Due to the stronger dissipation introduced by the source term, all the moments admit overall smaller scales, which can be observed from Figure \ref{fig:rte1d interm}. It is noticed that the algorithm generates accurate approximations to the first four moments, whereas the deviations in $m_4$ and $m_5$ appear significantly larger. We argue that the seemingly large errors are due to the excessively small moment magnitudes ($10^{-8}\sim 10^{-7}$), in which case the accumulated discretization errors from lower moments recovery might become dominating. In fact, the history of errors (Figure \ref{fig:rte1d interm err}) indicates the stable convergence of the approximate moments toward the exact solutions. 

We also mention that as reported in \cite{huang2022machineI}, the closure relation is relatively easy to capture in the intermediate and optically thick regimes, and the model \eqref{eq:ml closure} can achieve a good accuracy by only using  $m_0$, $m_1$, and $m_2$, which are accurately recovered in the computations. Hence, despite the difficulty in recovering high-order moments with extremely small magnitudes, we can still expect our algorithm to be effective when combined with actual machine-learning closure models.

\subsubsection{Boltzmann equation for channel flow}
In addition to the one-dimensional radiation transport, we also test our algorithm on particle flow through a closed channel. The flow is governed by a reduced Boltzmann equation with one dimension in space and two dimensions in velocity:
\begin{equation}\label{eq:boltz}
\begin{split}
\frac{\partial f}{\partial t}+v_1\partial_x f =  \iint_{\R^2\times \mathbb{S}^1}B(\bv-\bv_*,\sigma)(f'f_*'-ff_*) \,d\bv_* \,d\sig, \quad x\in\R, \hspace{0.5em} \bv = (v_1, v_2)\in\R^2,
\end{split}
\end{equation}
where
\begin{equation*}
\begin{split}
&f = f(x,\bv,t), \quad f_* = f(x,\bv_*,t), \quad f' = f(x,\bv',t), \quad f_*' = f(x,\bv_*',t),\\
\\
& \qquad \qquad \bv' = \frac{\bv+\bv_*}{2}+\frac{|\bv-\bv_*|}{2}\sig, \quad \bv_*' = \frac{\bv+\bv_*}{2}-\frac{|\bv-\bv_*|}{2}\sig.
\end{split}
\end{equation*}
The collision kernel $B(\bv-\bv_*,\sig)$ measures the intensity of particle collisions. In the case of the inverse power law, $B$ can be separated as the kinetic part $\Phi$ and angular part $b$,
    \[
		B(\bv-\bv_*,\sigma) = b(\cos\theta) \Phi(|\bv-\bv_*|), \quad \text{with} \ \cos\theta=\sigma \cdot \frac{\bv-\bv_*}{|\bv-\bv_*|},
    \]
where kinetic collision part $\Phi(|\bv-\bv_*|)=|\bv-\bv_*|^{\gamma}$ includes hard potential $ (\gamma>0) $, Maxwellian molecule $ (\gamma =0) $ and soft potential $ (\gamma<0) $, and angular part $b(\cos\theta)$ is often regarded to satisfy the Grad's cutoff assumption, i.e., $\int_{\mathbb{S}^{2}} b(\cos\theta) \,d \sigma < \infty$. For more details on the collision kernel, see \cite{Villani02}. For convenience of discussion, we consider the Maxwellian molecule with constant collision kernel $B= 0.2$. The initial condition is set to:
\[
\begin{split}
&f_0(x,v_1,v_2) = \rho(x)f_1(v_1)f_2(v_2), \quad x\in [0,1],\hspace{0.5em} (v_1,v_2)\in\R^2,
\end{split}
\]
with
\[
\begin{split}
&\rho(x) = 1+0.2\sin(2\pi x),\\
& f_1(v_1) = \frac{1}{\sig_1\sqrt{2\pi}}\exp(-\frac{(v_1-1)^2}{2\sig_1^2}), \quad v_1\in \R, \hspace{0.5em} \sig_1 = 0.5,\\
& f_2(v_2) = \frac{1}{\sig_2\sqrt{2\pi}}\exp(-\frac{(v_2-0.1)^2}{2\sig_2^2}), \quad v_2\in \R, \hspace{0.5em} \sig_2 = 0.05.
\end{split}
\]
Periodic boundary conditions are applied at both endpoints of the space interval. To generate reference data, the transport part is solved with WENO5-SSPRK3, while the collision part is calculated with the fast spectral method \cite{mouhot2006fast, HQ2020, HQY2021}.

\begin{figure}[h!]
    \centering
    \includegraphics[scale = 0.51]{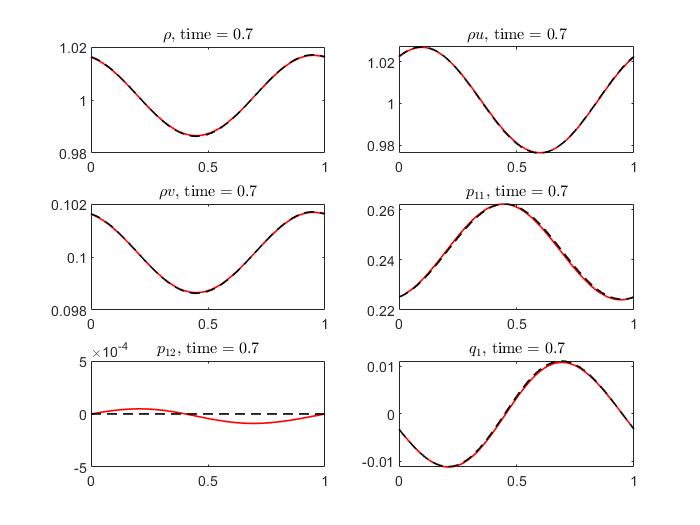}
    \caption{ Moment recovery for Boltzmann. Nudged solutions (red lines) vs. reference data (black dash lines). $t = 0.7$.}
    \label{fig:boltz1d}
\end{figure}

We define the moments of particle distribution $f$ up to the third order,
\[
\begin{split}
\rho(x,t) &:= \int_{\R^2} f(x,\bv,t) \,d\bv,\\ 
\rho\bu(x,t) &:= \int_{\R^2} \bv f(x,\bv,t)  \,d\bv, \qquad \bu = (u,v),\\
\pressure(x,t) &:=  \int_{\R^2}(\bv-\bu)\otimes(\bv-\bu)f(x,\bv,t) \,d\bv = \left[\begin{array}{cc}
    p_{11} & p_{12} \\
    p_{21} & p_{22}
\end{array}\right], \quad p_{12} = p_{21},\\
\rho e(x,t) &:= \int_{\R^2}\frac{|\bv-\bu|^2}{2}f(x,\bv,t) \,d\bv = \frac{1}{2}\trace(\pressure), \\
\bq(x,t) &:= \int_{\R^2}\frac{|\bv-\bu|^2}{2}(\bv-\bu)f(x,\bv,t) \,d\bv = [q_1,q_2]^\top.
\end{split}
\]
Integrating the equation \eqref{eq:boltz} against $1$, $\bv$, and $\displaystyle \frac{|\bv|^2}{2}$ yields the moment equations
\begin{equation}\label{eq:boltz moment 1d}
\left\{\begin{array}{l}
   \displaystyle \partial_t\rho+\partial_x(\rho u) = 0, \\
   \\
    \displaystyle \partial_t (\rho u) +\partial_x(\rho u^2)+\partial_x p_{11} = 0, \\
    \\
     \displaystyle \partial_t (\rho v) +\partial_x(\rho uv)+\partial_x p_{12} = 0,\\
    \\
    \displaystyle \partial_t(\rho e)+\partial_x(\rho eu)+\partial_x q_1+p_{11}\partial_x u+p_{12}\partial_x v= 0.
\end{array}
\right.
\end{equation}
Note that the internal energy is given by $\displaystyle \rho e = \frac{1}{2}(p_{11}+p_{22})$, while the evolution of $p_{22}$ is not described by the moment equations \eqref{eq:boltz moment 1d}. Nevertheless, due to the small variance of $f$ in $y$-velocity, the contribution of $p_{22}$ can be neglected in our test case, which yields $\displaystyle\rho e = \frac{1}{2}p_{11}$. Hence, the last equation in \eqref{eq:boltz moment 1d} reduces to 
\[
\partial_t p_{11}+\partial_x(p_{11} u)+2\partial_x q_1+2p_{11}\partial_x u+2p_{12}\partial_x v =0.
\]
Assume the complete observations of all moments are available at the boundaries. At the interior of the solution domain, we employ the partial observations with $\rho$, $\rho u$, and $\rho v$ over $N_{ob} = 60$ uniform observation grids. To recover the continuous data fields asymptotically in time, the solution procedures are briefly described as follows: at first, following the constructions of the nudged system \eqref{eq:nudge}, we reconstruct the approximate moments $\widehat{\rho}$, $\widehat{\rho u}$, $\widehat{\rho v}$, $\widehat{p}_{11}$, $\widehat{p}_{12}$ under $\mu = 10$. Then the approximate gradient of $q_1$ is obtained by
\[
\widetilde{\partial_x q}_1 = -\frac{1}{2}I_h(\partial_t \widehat{p}_{11}+\partial_x(\widehat{p}_{11}\widehat{u})+2\widehat{p}_{11}\partial_x\widehat{u}+2\widehat{p}_{12}\partial_x\widehat{v}),\quad \widehat{u}=\frac{\widehat{\rho u}}{\widehat{\rho}}, \hspace{0.5em} \widehat{v}= \frac{\widehat{\rho v}}{\widehat{\rho}}.
\]
Subsequently, the reconstruction of the third-order moment, $\widehat{q_1}$, is given by
\[
\widehat{q}_1 = q_1(0)+\int^x_0(\widetilde{\partial_x q}_1+c_q) \,dx, \quad c_q = -\int^1_0 \widetilde{\partial_x q_1} \,dx+ q_1(1)- q_1(0).
\]
The nudged equations are initialized with
\[
\widehat{\rho}_0(x) = 1, \quad \widehat{u}_0(x) = 1, \quad \widehat{v}_0(x) = 0.1, \quad \widehat{\pressure}_0(x) = \left[\begin{array}{cc}
    0 & 0 \\
    0 & 0
\end{array}\right] ,
\]
and discretized with WENO5-SSPRK3 under $N = 300$ space mesh points. Figure \ref{fig:boltz1d} displays the computational results at $t = 0.7$. The data of $\rho$, $\rho u$, $\rho v$, $p_{11}$ and $q_1$, which have major contributions to the system, are recovered with satisfactory accuracy. The seemingly large error of $\widehat{p}_{12}$ is, again, due to the extremely small magnitude of the exact $p_{12}$ ($\approx 10^{-9}$), hence the influence on the overall dynamics is minimal.

\section{Extension to higher dimension}
\label{sec:extension}

In the preceding discussions, we have demonstrated the efficacy of our moment recovery algorithm in one-dimensional space. Essentially, leveraging the sparsely observed data of lower-order moments enables the retrieval of gradients (variations) of higher-order moments at the interior of the domain. In addition, when complemented with sufficient ground truth data at the boundaries, the integration of the moment equations facilitates the recovery of all relevant moments.

However, extending the aforementioned procedure to recover higher-order moments in multiple dimensions is highly nontrivial. This complexity arises because the divergence term now involves a mixture of directional derivatives, which can lead to non-unique recovery of the quantity before differentiation. To illustrate this challenge and suggest a possible approach for addressing it, we briefly discuss the two-dimensional Eulerian system:
\[
\begin{split}
\left\{\begin{array}{l}
\displaystyle  \partial_t \rho + \nabla\cdot (\rho\bu) = 0 \,,\\
\\
\displaystyle \partial_t \rho\bu + \nabla\cdot(\rho\bu\otimes\bu)+\nabla\cdot\pressure = 0\,, \\
\\
\partial_t \rho e+\nabla\cdot(\rho e \bu)+\nabla\cdot \bq+\trace(\pressure\nabla\bu) = 0\,,
\end{array}\right. \quad \bx \in \Omega\in\R^2, \hspace{0.5em} t\geq 0.
\end{split}
\]
Assume the exact data of all solution variables are accessible at the boundary $\partial\Omega$, whereas the (sparse) interior observations are only available for $\rho$ and $\rho\bu$. We aim to recover the data of the pressure tensor, $\pressure\in\R^{2\times2}$, and subsequently the divergence of heat flux, $\nabla\cdot \bq$, using the continuous data assimilation. By evolving the corresponding nudged system along the lines of \eqref{eq:nudge}, one recovers the \emph{divergence} of the pressure, $g = -\nabla\cdot\pressure$, which contains the `mixture' of $x-$ and $y-$ derivatives. Hence, the recovery of $\pressure$ requires solving the divergence equations:
\begin{equation*}\label{eq: pressure div eqn}
\left\{\begin{array}{l}
     -(\partial_x p_{11}+\partial_y p_{21}) = g_1  \\
     \\
     -(\partial_x p_{12}+\partial_y p_{22}) = g_2
\end{array}\right., \quad p_{12} = p_{21},
\end{equation*}
which, in general, are under-determined. 

Indeed, denote $\pressure = [\pressure_1, \pressure_2]$ with $\pressure_j = [p_{1j},p_{2j}]^\top$, $j = 1, 2$, the pressure term can be represented via the Helmholtz decomposition:
\begin{equation}\label{eq:helmholtz}
\pressure_j  = \nabla\phi_j+\bR_j, \quad \bR_j = [R_{1j}, R_{2j}]^\top\,,
\end{equation}
where
\begin{subequations}
\begin{equation}\label{eq: phi poisson eqn}
\left\{\begin{array}{ll}
     -\Delta \phi_j = g_j,& \quad \bx\in \Omega\,,  \\
     \\
     \displaystyle
     \phi_j = 0,&\quad \bx \in \partial\Omega\,;
\end{array}\right.
\end{equation}
and
\begin{equation}\label{eq: R div eqn}
    \left\{\begin{array}{ll}
     \nabla\cdot \bR_j = 0,& \quad \bx\in \Omega\,,  \\
     \\
     \displaystyle
     \bR_j\cdot\bn = \pressure_j\cdot\bn-\nabla\phi_j\cdot\bn, & \quad \bx \in \partial\Omega\,.
\end{array}\right.
\end{equation}
\end{subequations}
Here $\bn = [n_1, n_2]^\top$ is the unit outer normal vector at the boundary. Moreover, the symmetry condition $p_{12} = p_{21}$ yields
\begin{equation}\label{eq: symm condition}
   R_{12} =R_{21}+\partial_y \phi_1- \partial_x\phi_2.
\end{equation}
Note that the Poisson equations \eqref{eq: phi poisson eqn} for the potentials $\phi_1$ and $\phi_2$ are well-posed. Therefore, the non-uniqueness of the pressure arises from the non-uniqueness of the solenoidal vector fields $\bR_j$. 

\begin{figure}
    \centering
    \begin{subfigure}{0.35\textwidth}
    \includegraphics[scale = 0.23]{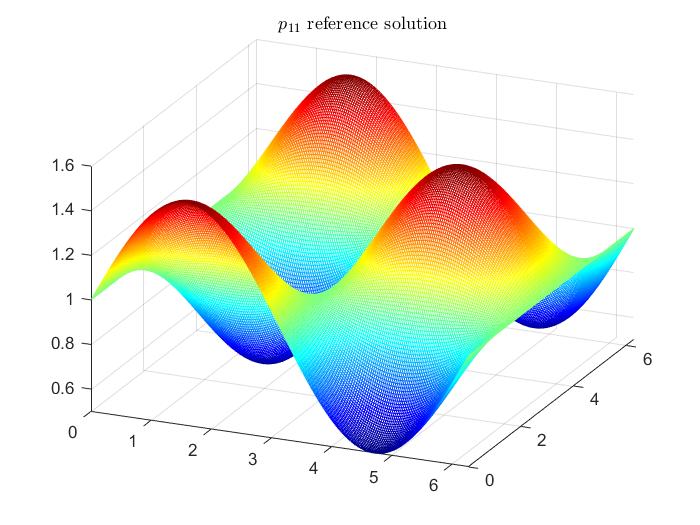}
    \subcaption{$p_{11}$}
    \end{subfigure}
    \quad
    \begin{subfigure}{0.35\textwidth}
    \includegraphics[scale = 0.23]{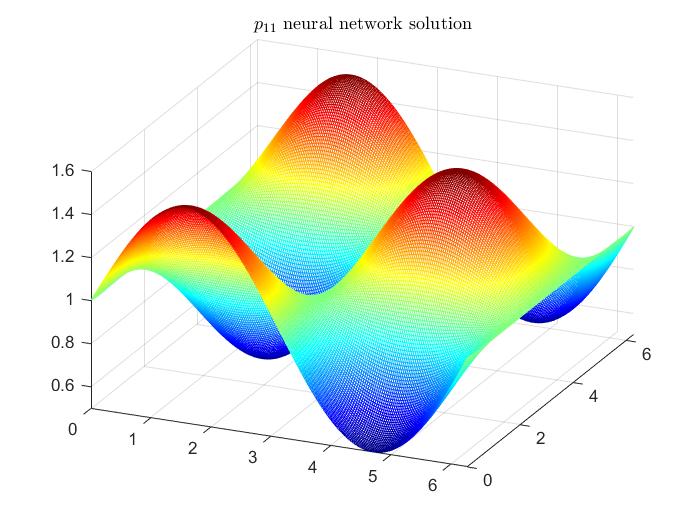}
    \subcaption{$p^{NN}_{11}$}
    \end{subfigure}
    \\
    \begin{subfigure}{0.35\textwidth}
    \includegraphics[scale = 0.23]{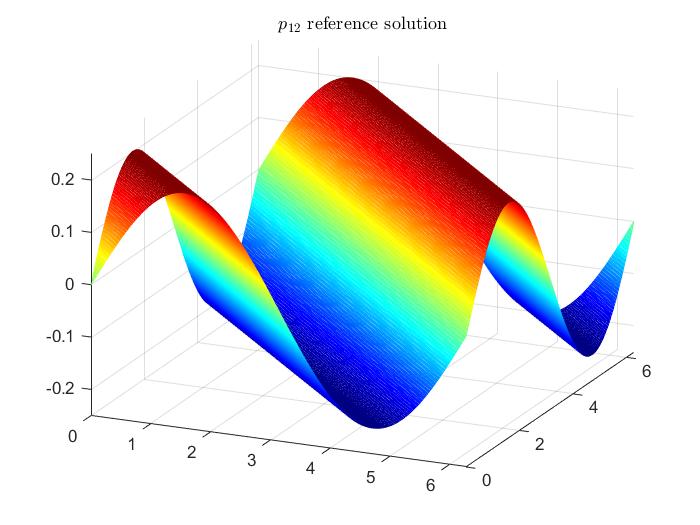}
    \subcaption{$p_{12}$}
    \end{subfigure}
    \quad
    \begin{subfigure}{0.35\textwidth}
    \includegraphics[scale = 0.23]{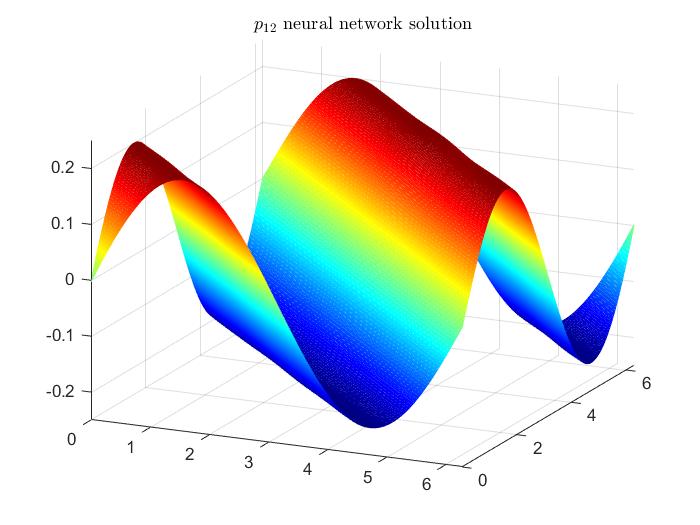}
    \subcaption{$p^{NN}_{12}$}
    \end{subfigure}
    \\
    \begin{subfigure}{0.35\textwidth}
    \includegraphics[scale = 0.23]{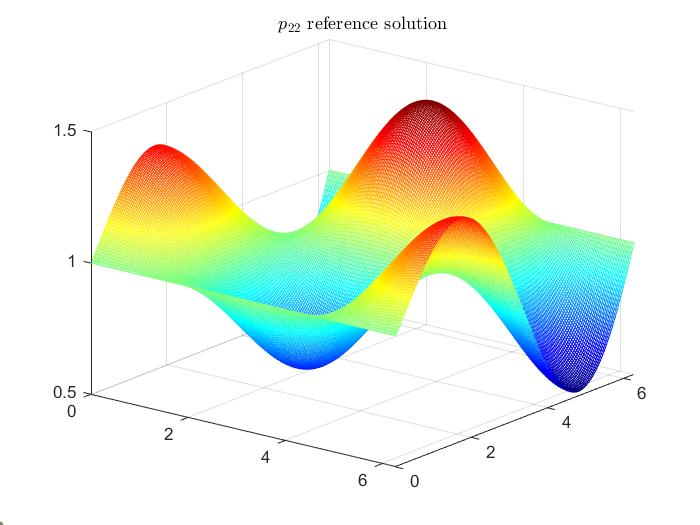}
    \subcaption{$p_{22}$}
    \end{subfigure}
    \quad
    \begin{subfigure}{0.35\textwidth}
    \includegraphics[scale = 0.23]{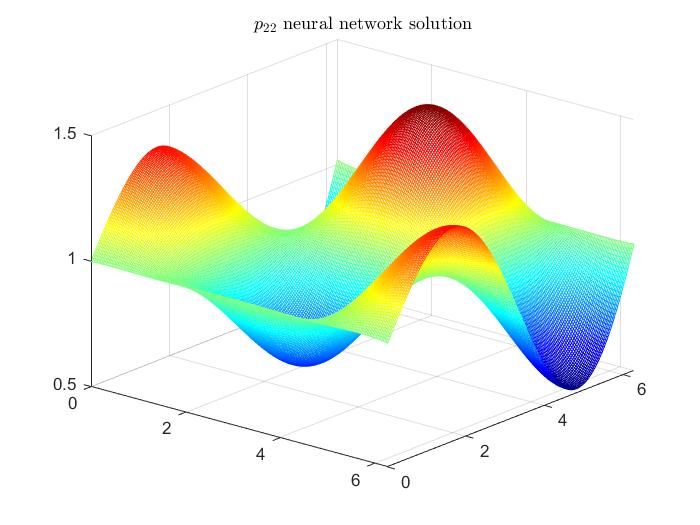}
    \subcaption{$p^{NN}_{22}$}
    \end{subfigure}
    \caption{Pressure recovery from divergence. Neural network solutions vs. reference solutions.}
    \label{fig:pressure recovery from div}
\end{figure}

For practical applications, a possible option is to generate approximate solutions with the physics-informed neural networks (PINNs) \cite{raissi2019physics}. Due to the simple linear structure of the divergence equations, the PINNs can be expected to work effectively. As an example, consider the exact pressure field $\pressure$ with
\[
\left\{\begin{array}{l}
     p_{11} = 1+0.5\sin(x)\cos(y)\,,  \\
     \\
     p_{12} = p_{21} =  0.2\sin(x+y)\,,\\
     \\
     p_{22} = 1+0.3e^{-0.1y}\cos(x)\sin(y)\,,
\end{array}\right. \quad (x,y)\in\Omega = [0,2\pi]^2.
\]
Based on the divergence at the interior, $\nabla\cdot\pressure_1 = \partial_x p_{11}+\partial_y p_{12}$ and $\nabla\cdot\pressure_2 = \partial_x p_{12}+\partial_y p_{22}$, and the ground truth boundary data, $\pressure|_{\partial \Omega}$, we reconstruct the pressure using the multi-outputs neural network,
\[
\network(x,y;\theta) = [p^{NN}_{11},p^{NN}_{12},p^{NN}_{22}]^\top,
\]
where $\network(\cdot;\theta)$ is a fully connected network with 3 hidden layers parameterized by $\th$. Each hidden layer has 20 neurons and the output is activated by the hyperbolic tangent function. The parameters of the neural network are optimized by minimizing the mean square loss,
\[
\begin{split}
MSE(\theta) &= \frac{1}{N}\sum^{N}_{j=1}|\partial_x p^{NN}_{11}(x_j,y_j)+\partial_y p^{NN}_{12}(x_j,y_j)-\nabla\cdot\pressure_1(x_j,y_j)|^2\\
&+\frac{1}{N}\sum^{N}_{j=1}|\partial_x p^{NN}_{12}(x_j,y_j)+\partial_y p^{NN}_{22}(x_j,y_j)-\nabla\cdot\pressure_2(x_j,y_j)|^2\\
&+\frac{1}{N_b}\sum^{2}_{m,n=1}\sum^{N_b}_{j=1}|p^{NN}_{mn}(x^b_j,y^b_j)-p_{mn}(x^b_j,y^b_j)|^2,
\end{split}
\]
where $\{(x_j,y_j)\}^N_{j=1}$ and $\{(x^b_j,y^b_j)\}^{N_b}_{j=1}$ are training points at the interior and the boundary. In our computation, we apply $30\times30$ uniformly distributed training points at the interior. At each of the four boundaries, $30$ training points are employed to enforce the boundary conditions. We train the network using LBFGS in Pytorch until the relative change in the loss is smaller than $10^{-4}$. As observed from the results displayed in Figure \ref{fig:pressure recovery from div}, the neural network solutions attain a reasonably good accuracy. The relatively errors are $\displaystyle \frac{||p^{NN}_{11}-p_{11}||_1}{||p_{11}||_1} = 1.54\times10^{-3}$, $\displaystyle \frac{||p^{NN}_{12}-p_{12}||_1}{||p_{12}||_1} = 7.05\times10^{-3}$, $\displaystyle \frac{||p^{NN}_{11}-p_{11}||_1}{||p_{11}||_1} = 1.51\times10^{-3}$. 

In general, to resolve the potential non-uniqueness issue, a regularization term can be added that favors the solution with special property. For instance, the solution with the smallest oscillation is favored by adding a penalty term $\displaystyle \lam||\nabla\pressure||^2_2$, where $\lam>0$ is a shrinking regularization parameter; by introducing the penalty $\displaystyle \lam||\bR||^2_2$, where $\bR = [\bR_1;\bR_2]$ are the divergence-free vector fields in the Helmholtz decomposition \eqref{eq:helmholtz}, one minimizes the projection of pressure onto $\image(\nabla)^\perp = \kernel(\nabla\cdot) = \{\mathbb{V}\in\R^{2\times2}:\ \nabla\cdot\mathbb{V} = 0, \ \mathbb{V} \cdot \bn|_{\partial\Omega} = \pressure \cdot \bn|_{\Omega}\}$.

Although using the PINNs provides a possible solution approach, one should be aware of the caveat in terms of the expensive computational cost. Further investigations are needed for more efficient data reconstruction from divergence information.

\bibliographystyle{amsplain}
\bibliography{Ref.bib}

\providecommand{\bysame}{\leavevmode\hbox to3em{\hrulefill}\thinspace}
\providecommand{\MR}{\relax\ifhmode\unskip\space\fi MR }
\providecommand{\MRhref}[2]{%
  \href{http://www.ams.org/mathscinet-getitem?mr=#1}{#2}
}
\providecommand{\href}[2]{#2}
\begin{thebibliography}{10}

\bibitem{albanez2018continuous}
D{\'e}bora~AF Albanez and Maicon~J Benvenutti, \emph{Continuous data assimilation algorithm for simplified bardina model}, Evolution Equations and Control Theory \textbf{7} (2018), no.~1, 33--52.

\bibitem{altaf2017downscaling}
M.~U. Altaf, E.~S. Titi, T.~Gebrael, O.~M. Knio, L.~Zhao, M.~F. McCabe, and I.~Hoteit, \emph{Downscaling the 2d b{\'e}nard convection equations using continuous data assimilation}, Computational Geosciences \textbf{21} (2017), 393--410.

\bibitem{azouani2014continuous}
Abderrahim Azouani, Eric Olson, and Edriss~S Titi, \emph{Continuous data assimilation using general interpolant observables}, Journal of Nonlinear Science \textbf{24} (2014), 277--304.

\bibitem{Bal_2009}
Guillaume Bal, \emph{Inverse transport theory and applications}, Inverse Problems \textbf{25} (2009), no.~5, 053001.

\bibitem{bardos1993fluid}
Claude Bardos, Fran{\c{c}}ois Golse, and C~David Levermore, \emph{Fluid dynamic limits of kinetic equations ii convergence proofs for the boltzmann equation}, Communications on pure and applied mathematics \textbf{46} (1993), no.~5, 667--753.

\bibitem{bhatnagar1954model}
Prabhu~Lal Bhatnagar, Eugene~P Gross, and Max Krook, \emph{A model for collision processes in gases. i. small amplitude processes in charged and neutral one-component systems}, Physical review \textbf{94} (1954), no.~3, 511.

\bibitem{cai2013globally}
Zhenning Cai, Yuwei Fan, and Ruo Li, \emph{Globally hyperbolic regularization of grad's moment system in one-dimensional space}, Communications in Mathematical Sciences \textbf{11} (2013), no.~2, 547--571.

\bibitem{cai2010numerical}
Zhenning Cai and Ruo Li, \emph{Numerical regularized moment method of arbitrary order for boltzmann-bgk equation}, SIAM Journal on Scientific Computing \textbf{32} (2010), no.~5, 2875--2907.

\bibitem{cai2012numerical}
Zhenning Cai, Ruo Li, and Yanli Wang, \emph{Numerical regularized moment method for high mach number flow}, Communications in Computational Physics \textbf{11} (2012), no.~5, 1415--1438.

\bibitem{carlson2020parameter}
Elizabeth Carlson, Joshua Hudson, and Adam Larios, \emph{Parameter recovery for the 2 dimensional navier--stokes equations via continuous data assimilation}, SIAM Journal on Scientific Computing \textbf{42} (2020), no.~1, A250--A270.

\bibitem{carlson2021dynamically}
Elizabeth Carlson, Joshua Hudson, Adam Larios, Vincent~R. Martinez, Eunice Ng, and Jared~P. Whitehead, \emph{Dynamically learning the parameters of a chaotic system using partial observations}, Discrete and Continuous Dynamical Systems. Series A \textbf{42} (2022), no.~8, 3809--3839. \MR{4447559}

\bibitem{carpenter1999improved}
James Carpenter, Peter Clifford, and Paul Fearnhead, \emph{Improved particle filter for nonlinear problems}, IEE Proceedings-Radar, Sonar and Navigation \textbf{146} (1999), no.~1, 2--7.

\bibitem{clark2018inferring}
Patricio Clark Di~Leoni, Andrea Mazzino, and Luca Biferale, \emph{Inferring flow parameters and turbulent configuration with physics-informed data assimilation and spectral nudging}, Physical Review Fluids \textbf{3} (2018), no.~10, 104604.

\bibitem{farhat2020data}
Aseel Farhat, Nathan~E Glatt-Holtz, Vincent~R Martinez, Shane~A McQuarrie, and Jared~P Whitehead, \emph{Data assimilation in large prandtl rayleigh--benard convection from thermal measurements}, SIAM Journal on Applied Dynamical Systems \textbf{19} (2020), no.~1, 510--540.

\bibitem{farhat2015continuous}
Aseel Farhat, Michael~S Jolly, and Edriss~S Titi, \emph{Continuous data assimilation for the 2d b{\'e}nard convection through velocity measurements alone}, Physica D: Nonlinear Phenomena \textbf{303} (2015), 59--66.

\bibitem{farhat2024identifying}
Aseel Farhat, Adam Larios, Vincent~R Martinez, and Jared~P Whitehead, \emph{Identifying the body force from partial observations of a two-dimensional incompressible velocity field}, Physical Review Fluids \textbf{9} (2024), no.~5, 054602.

\bibitem{farhat2016data}
Aseel Farhat, Evelyn Lunasin, and Edriss~S Titi, \emph{Data assimilation algorithm for 3d b{\'e}nard convection in porous media employing only temperature measurements}, Journal of Mathematical Analysis and Applications \textbf{438} (2016), no.~1, 492--506.

\bibitem{fischer2005overview}
Claude Fischer, Thibaut Montmerle, Lo{\"\i}k Berre, Ludovic Auger, and Simona~Ecaterina {\c{S}}tef{\u{a}}nescu, \emph{An overview of the variational assimilation in the aladin/france numerical weather-prediction system}, Quarterly Journal of the Royal Meteorological Society: A journal of the atmospheric sciences, applied meteorology and physical oceanography \textbf{131} (2005), no.~613, 3477--3492.

\bibitem{gesho2016computational}
Masakazu Gesho, Eric Olson, and Edriss~S Titi, \emph{A computational study of a data assimilation algorithm for the two-dimensional navier-stokes equations}, Communications in Computational Physics \textbf{19} (2016), no.~4, 1094--1110.

\bibitem{gottlieb2001strong}
Sigal Gottlieb, Chi-Wang Shu, and Eitan Tadmor, \emph{Strong stability-preserving high-order time discretization methods}, SIAM review \textbf{43} (2001), no.~1, 89--112.

\bibitem{HQ2020}
Jingwei Hu and Kunlun Qi, \emph{A fast {F}ourier spectral method for the homogeneous {B}oltzmann equation with non-cutoff collision kernels}, J. Comput. Phys. \textbf{423} (2020), 109806, 21. \MR{4156938}

\bibitem{HQY2021}
Jingwei Hu, Kunlun Qi, and Tong Yang, \emph{A new stability and convergence proof of the {F}ourier-{G}alerkin spectral method for the spatially homogeneous {B}oltzmann equation}, SIAM J. Numer. Anal. \textbf{59} (2021), no.~2, 613--633. \MR{4226997}

\bibitem{huang2022machineI}
Juntao Huang, Yingda Cheng, Andrew~J Christlieb, and Luke~F Roberts, \emph{Machine learning moment closure models for the radiative transfer equation i: directly learning a gradient based closure}, Journal of Computational Physics \textbf{453} (2022), 110941.

\bibitem{huang2023machineIII}
\bysame, \emph{Machine learning moment closure models for the radiative transfer equation iii: enforcing hyperbolicity and physical characteristic speeds}, Journal of Scientific Computing \textbf{94} (2023), no.~1, 7.

\bibitem{huang2023machineII}
Juntao Huang, Yingda Cheng, Andrew~J Christlieb, Luke~F Roberts, and Wen-An Yong, \emph{Machine learning moment closure models for the radiative transfer equation ii: Enforcing global hyperbolicity in gradient-based closures}, Multiscale Modeling \& Simulation \textbf{21} (2023), no.~2, 489--512.

\bibitem{jolly2019continuous}
Michael~S Jolly, Vincent~R Martinez, Eric~J Olson, and Edriss~S Titi, \emph{Continuous data assimilation with blurred-in-time measurements of the surface quasi-geostrophic equation}, Chinese Annals of Mathematics, Series B \textbf{40} (2019), no.~5, 721--764.

\bibitem{jolly2017data}
Michael~S Jolly, Vincent~R Martinez, and Edriss~S Titi, \emph{A data assimilation algorithm for the subcritical surface quasi-geostrophic equation}, Advanced Nonlinear Studies \textbf{17} (2017), no.~1, 167--192.

\bibitem{jolly2015determining}
Michael~S Jolly, Tural Sadigov, and Edriss~S Titi, \emph{A determining form for the damped driven nonlinear schr{\"o}dinger equation—fourier modes case}, Journal of Differential Equations \textbf{258} (2015), no.~8, 2711--2744.

\bibitem{jolly2017determining}
\bysame, \emph{Determining form and data assimilation algorithm for weakly damped and driven korteweg--de vries equation—fourier modes case}, Nonlinear Analysis: Real World Applications \textbf{36} (2017), 287--317.

\bibitem{kalman1960new}
Rudolph~Emil Kalman, \emph{A new approach to linear filtering and prediction problems},  (1960).

\bibitem{kitagawa1996monte}
Genshiro Kitagawa, \emph{Monte carlo filter and smoother for non-gaussian nonlinear state space models}, Journal of computational and graphical statistics \textbf{5} (1996), no.~1, 1--25.

\bibitem{le1986variational}
Fran{\c{c}}ois-Xavier Le~Dimet and Olivier Talagrand, \emph{Variational algorithms for analysis and assimilation of meteorological observations: theoretical aspects}, Tellus A: Dynamic Meteorology and Oceanography \textbf{38} (1986), no.~2, 97--110.

\bibitem{martinez2022convergence}
Vincent~R Martinez, \emph{Convergence analysis of a viscosity parameter recovery algorithm for the 2d navier--stokes equations}, Nonlinearity \textbf{35} (2022), no.~5, 2241.

\bibitem{martinez2022reconstruction}
Vincent~R. Martinez, \emph{On the reconstruction of unknown driving forces from low-mode observations in the 2d navier–stokes equations}, Proceedings of the Royal Society of Edinburgh: Section A Mathematics (2024), 1–24.

\bibitem{mouhot2006fast}
Cl{\'e}ment Mouhot and Lorenzo Pareschi, \emph{Fast algorithms for computing the boltzmann collision operator}, Mathematics of computation \textbf{75} (2006), no.~256, 1833--1852.

\bibitem{pachev2022concurrent}
Benjamin Pachev, Jared~P Whitehead, and Shane~A McQuarrie, \emph{Concurrent multiparameter learning demonstrated on the kuramoto--sivashinsky equation}, SIAM Journal on Scientific Computing \textbf{44} (2022), no.~5, A2974--A2990.

\bibitem{pei2018continuous}
Yuan Pei, \emph{Continuous data assimilation for the 3d primitive equations of the ocean}, arXiv preprint arXiv:1805.06007 (2018), 1--25.

\bibitem{raissi2019physics}
Maziar Raissi, Paris Perdikaris, and George~E Karniadakis, \emph{Physics-informed neural networks: A deep learning framework for solving forward and inverse problems involving nonlinear partial differential equations}, Journal of Computational physics \textbf{378} (2019), 686--707.

\bibitem{Villani02}
C\'{e}dric Villani, \emph{A review of mathematical topics in collisional kinetic theory}, Handbook of Mathematical Fluid Mechanics (S.~Friedlander and D.~Serre, eds.), vol.~I, North-Holland, 2002, pp.~71--305.

\bibitem{welch1995introduction}
Greg Welch and Gary Bishop, \emph{An introduction to the kalman filter},  (1995).

\bibitem{vzupanski1995four}
Du{\v{s}}anka {\v{Z}}upanski and Fedor Mesinger, \emph{Four-dimensional variational assimilation of precipitation data}, Monthly Weather Review \textbf{123} (1995), no.~4, 1112--1127.

\end{thebibliography}






\end{document}